\documentclass[preprint,12pt,authoryear]{elsarticle}

\biboptions{authoryear,round}

\usepackage[a4paper, margin=2.5cm]{geometry}
\usepackage{microtype}
\usepackage{amsmath}
\usepackage{amssymb}
\usepackage{amsthm}
\usepackage{bm}
\usepackage{cases}

\usepackage{graphicx}
\usepackage{xcolor}
\definecolor{linkblue}{RGB}{0,70,140}
\usepackage{booktabs}

\usepackage{float}
\usepackage{algorithm}
\usepackage{algpseudocode}

\usepackage{etoolbox}
\newtheoremstyle{boldremark}
{3pt}
{3pt}
{\normalfont}
{}
{\bfseries}
{.}
{ }
{}

\theoremstyle{boldremark}
\newtheorem{remark}{Remark}
\AtEndEnvironment{remark}{\unskip\nobreak\hfill\(\square\)\par}

\usepackage{xurl}
\usepackage[
colorlinks=true,
linkcolor=linkblue,
citecolor=linkblue,
urlcolor=linkblue,
breaklinks=true
]{hyperref}

\usepackage[nameinlink,capitalise,noabbrev]{cleveref}

\crefname{figure}{Fig.}{Figs.}
\Crefname{figure}{Fig.}{Figs.}

\crefname{table}{Table}{Tables}
\Crefname{table}{Table}{Tables}

\crefname{equation}{Eq.}{Eqs.}
\Crefname{equation}{Eq.}{Eqs.}

\crefname{section}{Section}{Sections}
\Crefname{section}{Section}{Sections}

\crefname{subsection}{Section}{Sections}
\Crefname{subsection}{Section}{Sections}

\crefname{subsubsection}{Section}{Sections}
\Crefname{subsubsection}{Section}{Sections}

\crefname{remark}{Remark}{Remarks}
\Crefname{remark}{Remark}{Remarks}

\crefname{algorithm}{Algorithm}{Algorithms}
\Crefname{algorithm}{Algorithm}{Algorithms}

\begin{document}
	
	\begin{frontmatter}
		
		\title{A shifted energy barrier approach for phase-field modeling of tensile-dominated brittle fracture} 
		
		\author[inst1,inst2]{Yaode Yin}
		\author[inst2]{Luigi Greco}
		\author[inst1]{Hongjun Yu\corref{cor1}}
		\ead{yuhongjun@hit.edu.cn}
		\author[inst2]{Simone Morganti}
		\cortext[cor1]{Corresponding author.}
		\affiliation[inst1]{
			organization={Department of Astronautic Science and Mechanics},
			addressline={Harbin Institute of Technology},
			city={Harbin},
			country={China}
		}
		
		\affiliation[inst2]{
			organization={Department of Civil Engineering and Architecture},
			addressline={University of Pavia},
			city={Pavia},
			country={Italy}
		}
		
		\begin{abstract}
			The classical \(\mathrm{AT}_1\) phase-field model contains an intrinsic energy barrier for crack nucleation, which makes the predicted strength depend on the fracture toughness and the regularization length. For tensile-dominated brittle fracture, this barrier is shifted by mapping the Rankine criterion, evaluated on the effective stress, onto a state-dependent active-energy threshold. The prescribed tensile strength then controls crack nucleation, while the \(\mathrm{AT}_1\) crack-density functional, stiffness degradation, and degraded stress response remain unchanged. Since the threshold depends on the current stress state, the field equations are derived from a restricted variational principle. A microforce formulation identifies the barrier shift as a dissipative resistance and provides the corresponding lower bound on the regularization length. In one-dimensional tension, closed-form solutions recover the prescribed peak strength and give a cosine-type localization profile that approaches the classical \(\mathrm{AT}_1\) profile as the shift vanishes. Numerical examples show that, within the admissible range, the nucleation load is nearly insensitive to the regularization length and the predicted multiaxial nucleation states follow the Rankine envelope. Under overall compression, crack nucleation remains associated with local tensile stress concentrations. The formulation also captures the transition from strength-controlled failure for small flaws to the LEFM limit for large cracks.
		\end{abstract}
		
		
	\end{frontmatter}

	\section{Introduction}
	\label{sec:intro}
	
	Fracture modeling has been developed through a continuous interaction between physical fracture concepts and numerical descriptions of evolving cracks. 
	The Griffith theory provided the energetic basis for brittle fracture by relating crack propagation to the balance between released elastic energy and the energy required to create new crack surfaces \citep{griffith1921vi}. 
	This theory established the critical energy release rate as the central quantity for crack growth. 
	At the same time, the onset of fracture in an initially intact body also involves local strength and the formation of a fracture process zone \citep{orowan1949fracture,Irwin1960,leguillon2002strength,cornetti2006finite}.  
	Cohesive-zone models introduced this aspect through traction--separation relations, thereby connecting crack opening, cohesive strength, softening behavior, and fracture energy \citep{Barenblatt1962}. 
	Their numerical implementation in cohesive surface formulations and element-based crack insertion methods provided important tools for simulating crack initiation and dynamic crack growth \citep{Xu1994,Camacho1996}. 
	Other sharp-crack approaches, such as extended finite element methods, represented cracks through enriched displacement approximations and allowed crack growth without conforming remeshing \citep{Moes1999}. 
	These developments form an important background for modern computational fracture mechanics and highlight the different ways in which crack surfaces, process zones, and evolving discontinuities can be represented.
	
	A different and now widely used route is provided by the variational phase-field fracture method. 
	Instead of treating the crack as an explicit sharp discontinuity, phase-field models approximate the crack set by a smooth scalar field and describe fracture evolution through a regularized energy functional. 
	The theoretical foundation of this approach is closely related to the variational reformulation of Griffith fracture by \citet{Francfort1998}, in which brittle fracture is written as an energy minimization problem involving bulk elastic energy and crack surface energy. 
	\citet{Bourdin2000} introduced a diffusive approximation of the crack set that made this variational formulation suitable for finite element implementation, and \citet{Bourdin2007} further developed its numerical realization for quasi-static brittle fracture. 
	Existence and convergence results for quasi-static evolutions provided an important mathematical basis for the variational setting \citep{Francfort2003}. 
	The connection between phase-field fracture and gradient-damage models was clarified by \citet{Pham2011a}, who showed how regularized damage functionals can approximate brittle fracture, and by \citet{Pham2011b}, who analyzed the stability and uniqueness of homogeneous responses in phase-field models. 
	The tension--compression asymmetry of the elastic energy and the associated numerical solution strategies were further shaped by the works of \citet{Amor2009}, \citet{Miehe2010a}, and \citet{Miehe2010b}. 
	These developments established the classical phase-field framework and its commonly used \(\mathrm{AT}_1\) and \(\mathrm{AT}_2\) variants, as summarized in later reviews and comparative studies \citep{Ambati2015review,DeBorst2016,wu2020phase,greco2026fourth}.
	
	Within the variational phase-field framework, the growth of a pre-existing crack is naturally connected to the Griffith energy balance. 
	The regularized crack-surface functional introduced by \citet{Bourdin2000} provides a diffuse approximation of the sharp crack surface, and the resulting formulation is consistent with the Francfort--Marigo variational theory of brittle fracture in the sharp-crack limit \citep{Francfort1998,Francfort2003,Bourdin2007,negri2020gamma}. 
	For problems with existing cracks or sufficiently strong stress concentrations, the critical energy release rate controls the energetic cost of creating additional crack surface. 
	This is the setting in which Griffith-type phase-field models have been particularly successful.
	Crack nucleation in an initially intact body involves an additional issue. 
	Before a macroscopic crack has formed, there is no prescribed crack front or crack extension on which a classical Griffith-type energy release rate can be directly evaluated. 
	In a regularized phase-field model, the onset of damage is instead associated with the stability of the homogeneous intact or weakly damaged state and with the subsequent localization of the phase-field variable \citep{Pham2011b,Tanne2018,Zolesi2024}. 
	The predicted nucleation load is therefore not determined by the critical energy release rate alone. 
	It also depends on the crack-density function, the degradation function, the adopted tension--compression split, the regularization length, and the geometry-induced stress concentration \citep{Pham2011a,Pham2011b,Tanne2018,DeLorenzis2022,Vicentini2024}.

	This kind of distinction is important because critical energy release rate and strength are usually regarded as different material properties. 
	The critical energy release rate characterizes the energy required for crack growth, whereas the strength characterizes the stress level associated with failure initiation in a nominally intact material \citep{LopezPamies2025}. 
	In the standard \(\mathrm{AT}_1\) formulation, damage does not start immediately from the undegraded state \citep{Pham2011a}, which makes the model attractive for studying delayed crack initiation \citep{DeLorenzis2022,Vicentini2024,Kumar2018}. 
	However, this elastic stage is produced by an intrinsic barrier of the regularized crack-surface functional. 
	For the usual crack density and quadratic degradation function, this barrier is controlled by elastic constants, critical energy release rate and the regularization length \citep{miehe2015phase}. 
	As a result, the apparent nucleation strength of the standard model is tied to the combination of elastic constants, critical energy release rate and the regularization length, rather than being an independently prescribed material parameter \citep{Tanne2018,Kristensen2020}.

	The observation of length-scale dependence has driven extensive research into the control of phase-field crack nucleation. 
	Early efforts focused on one-dimensional analytical solutions, homogeneous uniaxial tension responses, and tensile-dominated benchmark problems. 
	\citet{Tanne2018} systematically analyzed crack nucleation from homogeneous states and in structural geometries, showing that the predicted nucleation load depends on the regularization length, the geometry, and the severity of stress concentration. 
	To decouple the macroscopic material strength from the regularization length, \citet{lorentz2011convergence} formulated a gradient-damage approach that converges asymptotically to a Barenblatt cohesive-zone model. 
	Building on related asymptotic and cohesive interpretations, \citet{Wu2017} and \citet{Wu2018} developed length-scale-insensitive phase-field damage models by introducing suitable crack geometric functions and rational degradation functions. 
	\citet{Feng2021} further generalized cohesive phase-field modeling using the integral transform, allowing prescribed cohesive laws to be embedded into the phase-field fracture framework. 
	In these approaches the strength enters through one-dimensional analytical solutions, in which the nucleation stress is expressed in terms of the elastic constants, the fracture toughness, the regularization length and the constitutive functions. The calibration is therefore anchored to the uniaxial response.
	
	Strength control under complex multiaxial loading requires additional constitutive considerations. 
	Classical brittle phase-field formulations commonly rely on tension--compression energy decompositions, such as the volumetric--deviatoric split \citep{Amor2009} and the spectral split \citep{Miehe2010a,Miehe2010b}, to prevent damage evolution under compression. 
	However, \citet{DeLorenzis2022} and \citet{Vicentini2024} demonstrated that standard energy decompositions also act as implicit constitutive assumptions on the multiaxial failure envelope once a model is calibrated to tensile strength. 
	Furthermore, the stability analyses of \citet{Zolesi2024} showed that the onset of localization under multiaxial loading depends on how different stress components are allowed to drive damage through the active part of the elastic energy. 
	Among the recently developed energy decompositions, the star-convex decomposition \citep{Vicentini2024} and the directional energy decomposition \citep{feng2023unified} are two representative examples. 
	The former extends the volumetric--deviatoric split of \citet{Amor2009} by introducing an additional parameter, thereby increasing the flexibility of the implied strength surface while retaining the standard variational structure. 
	The latter is based on structured-deformation theory \citep{freddi2010regularized,tadmor2011modeling} and determines the macroscopic elastic energy density through a homogenization procedure, thereby allowing the crack direction to be selected by energy minimization while the tensile-to-shear strength ratio influences the predicted fracture mode and the orientation of crack initiation.

	Another route is to incorporate an explicit strength criterion into the Griffith-type phase-field theory. 
	\citet{Kumar2018} and \citet{Kumar2020} emphasized that fracture nucleation requires material strength as an independent ingredient, distinct from the critical energy release rate. 
	In their formulation, the nucleation process is represented through an additional phase-field fracture driving force associated with the macroscopic effect of microscopic defects. 
	\citet{Larsen2024} later recast this strength-based idea in a variational form, showing how Griffith phase-field fracture with material strength can be formulated through separate minimization structures. 
	More recently, \citet{LopezPamies2025} sharpened this viewpoint by arguing that classical variational phase-field models require an independent strength input if they are to predict fracture nucleation in a physically meaningful way.

	A more recent class of formulations departs more fundamentally from the standard stiffness-degradation structure of \(\mathrm{AT}\)-type phase-field models. 
	Instead of applying the degradation function directly to the elastic energy, these models introduce the strength surface through a support-function-based energy construction and let the phase-field variable degrade the admissible strength domain. 
	\citet{Vicentini2025} proposed a variational cohesive phase-field model based on an eigenstrain formulation, in which an arbitrary convex strength surface can be prescribed independently of the regularization length and the cohesive response can be tuned flexibly. 
	\citet{Bourdin2025} developed a variational framework for fracture with arbitrary closed convex strength domains, where the material stiffness is not degraded; instead, the strength domain contracts as the phase-field variable evolves. 
	The key distinction of these formulations is therefore a shift from stiffness degradation to strength degradation, which allows strength criteria to enter the variational phase-field theory as constitutive input rather than emerging indirectly from a degraded elastic energy.
	
	Among the formulations reviewed above, those that introduce material strength as an explicit constitutive input under multiaxial conditions follow two main routes. In the external driving force approach \citep{Kumar2020,Larsen2024,LopezPamies2025}, the phase-field evolution equation is supplemented with an additional driving force constructed from the strength surface. In the strength-degradation approach \citep{Vicentini2025,Bourdin2025}, the stiffness-degradation structure is replaced and the phase-field variable contracts the admissible strength domain. Both routes provide the flexibility to incorporate general convex criteria. For applications focused on tensile-dominated brittle fracture, a complementary question remains open: can the macroscopic tensile strength be prescribed under multiaxial stress states through a direct mapping of the stress criterion to the damage-initiation threshold, while the \(\mathrm{AT}_1\) crack surface density, the stiffness-degradation structure and the standard degraded stress response are all kept unchanged?
	
	To address this, the present work proposes a Rankine-shifted \(\mathrm{AT}_1\) energy barrier. The core mechanism relies on a direct stress-to-energy mapping rather than a reconstruction of the global energy functional. Specifically, a Rankine criterion governed by the maximum principal effective stress is mapped to a local active energy threshold. By shifting the intrinsic \(\mathrm{AT}_1\) nucleation barrier with the strength-based threshold, damage initiation is governed directly by the Rankine criterion. Concurrently, the standard crack-surface density function and the conventional degraded elastic stress response remain completely intact. Because the activation threshold adapts to the pre-critical stress state, the governing equations are cast in a restricted variational principle \citep{rosen1954solution,rosen1954use}. 
	
	The remainder of the paper is organized as follows. 
	\cref{sec:classical} reviews the classical phase-field formulation and derives the intrinsic \(\mathrm{AT}_1\) barrier for damage initiation. 
	\cref{sec:shifted_theory} introduces the Rankine-equivalent threshold, the shifted damage-driving force, and the restricted variational formulation. 
	\cref{sec:thermodynamic_admissibility} discusses the freezing rule and the thermodynamic admissibility of the shifted damage evolution. 
	\cref{sec:1d_analysis} specializes the formulation to the uniaxial traction problem of a one-dimensional bar, in which the Rankine-equivalent threshold reduces to a constant. The resulting closed-form solutions isolate the effect of the barrier shift and provide analytical benchmarks for the numerical implementation.
	\cref{sec:numerical_implementation} presents the numerical solution strategy. 
	\cref{sec:examples} examines the dependence on the regularization length, the response under compression, the transition from strength-controlled to toughness-controlled failure, and multiaxial crack initiation. 
	Finally, \cref{sec:conclusions} summarizes the main findings and limitations of the formulation.

	\section{Classical phase-field fracture method and its intrinsic barrier}
	\label{sec:classical}
	
	\subsection{Variational approach to Griffith's fracture}
	\label{subsec:variational_griffith}
	
	Let $\Omega \subset \mathbb{R}^n$ ($n \in \{1,2,3\}$) be an open and bounded domain representing a solid body, with its external boundary denoted by $\partial\Omega$. An internal sharp crack is represented by a lower-dimensional set $\mathcal{C}$ embedded in the body. Following \citet{Francfort1998}, Griffith's brittle fracture theory can be formulated as an energy minimization problem.
	
	Neglecting body forces and surface tractions for simplicity, and considering a time-dependent Dirichlet boundary condition $\mathbf{u}(\mathbf{x},t) = \bar{\mathbf{u}}(\mathbf{x},t)$ on $\partial_D \Omega$, the total potential energy of the cracked body is written as
	\begin{equation}
		\label{eq:griffith}
		\mathcal{E}_{\mathcal{C}}\left(\mathbf{u},\mathcal{C}\right)
		=
		\int_{\Omega \setminus \mathcal{C}}
		\psi_e(\bm{\varepsilon}(\mathbf{u})) \, d\Omega
		+
		\int_{\mathcal{C}} G_c \, d\mathcal{S},
	\end{equation}
	where $G_c > 0$ is the critical energy release rate, and $\mathbf{u}$ is the kinematically admissible displacement field. The infinitesimal strain tensor is defined as $\bm{\varepsilon}(\mathbf{u}) = \frac{1}{2}(\nabla \mathbf{u} + (\nabla \mathbf{u})^T)$.
	
	For an isotropic linear elastic material, the elastic strain energy density is
	\begin{equation}
		\label{eq:elastic_energy}
		\psi_e(\bm{\varepsilon})
		=
		\frac{\kappa}{2}\text{tr}(\bm\varepsilon)^2
		+
		\mu \bm{\varepsilon}_{\mathrm{dev}}:\bm{\varepsilon}_{\mathrm{dev}},
	\end{equation}
	where $\kappa$ and $\mu$ are the bulk and shear moduli, respectively. The strain tensor is decomposed into spherical and deviatoric parts as
	\begin{equation}
		\bm\varepsilon
		=
		\bm\varepsilon_{\mathrm{sph}}
		+
		\bm\varepsilon_{\mathrm{dev}},
		\qquad
		\bm\varepsilon_{\mathrm{sph}}
		=
		\frac{\text{tr}(\bm\varepsilon)}{n}\mathbf{I},
		\qquad
		\bm\varepsilon_{\mathrm{dev}}
		=
		\bm\varepsilon
		-
		\frac{\text{tr}(\bm\varepsilon)}{n}\mathbf{I},
	\end{equation}
	where $\mathbf{I}$ is the second-order identity tensor. The Cauchy stress tensor follows from the elastic potential:
	\begin{equation}
		\bm{\sigma}_0(\bm{\varepsilon})
		=
		\frac{\partial\psi_e(\bm{\varepsilon})}{\partial\bm{\varepsilon}}
		=
		\kappa\text{tr}(\bm{\varepsilon})\mathbf{I}
		+
		2\mu\bm\varepsilon_{\mathrm{dev}}.
	\end{equation}
	
	Within this variational setting, the quasi-static evolution of $\mathbf{u}$ and $\mathcal{C}$ at time $t$ is obtained by minimizing \cref{eq:griffith} subject to the crack irreversibility constraint:
	\begin{equation}
		\label{eq:argloc}
		(\mathbf{u},\mathcal{C})
		=
		\arg \min_{\mathbf{u}, \mathcal{C}}
		\quad
		\mathcal{E}_\mathcal{C}(\mathbf{u},\mathcal{C}),
		\quad
		\text{s.t.}
		\quad
		\mathcal{C}(t) \supseteq \mathcal{C}(s)
		\quad
		\text{for } t > s.
	\end{equation}
	
	\begin{remark}
		The unilateral constraint $\mathcal{C}(t) \supseteq \mathcal{C}(s)$, often expressed in rate form as $\dot{\mathcal{C}} \geq 0$, enforces crack irreversibility and prevents crack healing during the quasi-static process.
	\end{remark}
	
	\subsection{Phase-field approximation and energy decomposition}
	\label{subsec:phase_field_approx}
	
	Tracking the discontinuous crack surface $\mathcal{C}$ in a complex domain is computationally demanding. The phase-field method avoids explicit crack tracking by introducing a continuous scalar damage variable $d(\mathbf{x},t)\in[0,1]$, where $d=0$ denotes the intact state and $d=1$ denotes the fully broken state. The sharp crack is regularized over a finite width controlled by the regularization length $l_0>0$.
	
	Among phase-field formulations, the \(\mathrm{AT}_1\) model is useful for studying crack nucleation because it provides an elastic limit before damage initiation. The Griffith surface energy in \cref{eq:griffith} is approximated by the \(\mathrm{AT}_1\) crack density functional \citep{Pham2011a,Tanne2018}:
	\begin{equation}
		\label{eq:at1_fracture_energy}
		\int_{\mathcal{C}} G_c \, d\mathcal{S}
		\approx
		\frac{3G_c}{8l_0}
		\int_{\Omega}
		\left(
		d
		+
		l_0^2 |\nabla d|^2
		\right)
		d\Omega.
	\end{equation}
	
	To reduce non-physical crack growth under compression, a tension--compression split of the elastic energy is commonly introduced. Following \citet{Amor2009}, the strain energy density is decomposed into an active part $\psi_e^+$ and a passive part $\psi_e^-$:
	\begin{equation}
		\label{eq:amor_split}
		\psi_e(\bm\varepsilon)
		=
		\psi_e^+(\bm\varepsilon)
		+
		\psi_e^-(\bm\varepsilon),
	\end{equation}
	with
	\begin{equation}
		\psi_e^+(\bm\varepsilon)
		=
		\frac{\kappa}{2}
		\langle \text{tr}(\bm\varepsilon) \rangle_+^2
		+
		\mu \bm{\varepsilon}_{\mathrm{dev}} : \bm{\varepsilon}_{\mathrm{dev}},
		\qquad
		\psi_e^-(\bm\varepsilon)
		=
		\frac{\kappa}{2}
		\langle \text{tr}(\bm\varepsilon) \rangle_-^2,
	\end{equation}
	where $\langle x \rangle_{\pm} = (x \pm |x|)/2$ are the Macaulay brackets. The corresponding stress contributions are denoted by $\bm\sigma_0^+ = \partial\psi_e^+ / \partial\bm\varepsilon$ and $\bm\sigma_0^- = \partial\psi_e^- / \partial\bm\varepsilon$. \citet{Vicentini2024} established that eliminating spurious residual shear stresses at full damage is essential to prevent unphysical traction transmission across separated crack faces, a physical requirement strictly satisfied by Amor's decomposition.
	
	Using the degradation function $g(d)=(1-d)^2$, the regularized \(\mathrm{AT}_1\) energy functional is
	\begin{equation}
		\label{eq:at1_total_energy}
		\mathcal{E}_{\mathrm{AT}_1}^*(\mathbf{u},d)
		=
		\int_{\Omega}
		\Bigl(
		g(d)\psi_e^+(\bm{\varepsilon}(\mathbf{u}))
		+
		\psi_e^-(\bm{\varepsilon}(\mathbf{u}))
		\Bigr)
		d\Omega
		+
		\frac{3G_c}{8l_0}
		\int_{\Omega}
		\left(
		d
		+
		l_0^2 |\nabla d|^2
		\right)
		d\Omega.
	\end{equation}
	
	For a quasi-static process at time \(t\in[0,T]\), the coupled state \((\mathbf{u}(\mathbf{x},t),d(\mathbf{x},t))\) is obtained from the constrained local minimization problem
	\begin{equation}
		\label{eq:local_min_continuous}
		(\mathbf{u}, d)
		=
		\arg \text{loc} \min_{(\mathbf{u}^*, d^*) \in \mathcal{U}_t \times \mathcal{D}_t}
		\mathcal{E}_{\mathrm{AT}_1}^*(\mathbf{u}^*, d^*),
	\end{equation}
	where
	\begin{equation}
		\label{eq:space_u_cont}
		\mathcal{U}_t
		=
		\left\{
		\mathbf{u}^* \in H^1(\Omega; \mathbb{R}^n)
		\mid
		\mathbf{u}^*(\mathbf{x}) = \bar{\mathbf{u}}(\mathbf{x}, t)
		\text{ on } \partial_D\Omega
		\right\},
	\end{equation}
	and
	\begin{equation}
		\label{eq:space_d_cont}
		\mathcal{D}_t
		=
		\left\{
		d^* \in H^1(\Omega; \mathbb{R})
		\mid
		\max_{s \in [0, t)} d(\mathbf{x}, s)
		\leq
		d^*(\mathbf{x})
		\leq
		1
		\text{ a.e. in } \Omega
		\right\}.
	\end{equation}
	Here, $\mathcal{D}_t$ enforces the irreversibility of the phase field by requiring the current damage field to be no smaller than its previous maximum value.
	
	\begin{remark}
		\label{rem:tension_driven}
		In \cref{eq:at1_total_energy}, only $\psi_e^+$ is degraded by $g(d)$, so 
		damage is driven by the active energy alone. Amor's split puts the whole 
		deviatoric energy into $\psi_e^+$. A shear-dominated state can then build up 
		enough active energy to nucleate damage even when no principal stress is 
		tensile. In the brittle solids we consider, crack onset is instead set by 
		the maximum tensile stress reaching the strength 
		\citep{leguillon2002strength, cornetti2006finite}, so this shear-driven 
		nucleation does not match the tensile-dominated failure we want. This is 
		what motivates the Rankine-based shift below, which changes the onset condition of damage, tying it to $\sigma_1 = \sigma_t$. We keep Amor's split otherwise 
		unchanged: $\psi_e^-$ is left undegraded, so closed crack faces still carry 
		compression and do not interpenetrate.
	\end{remark}
	
	\subsection{The intrinsic energy barrier for crack nucleation}
	\label{subsec:inherent_barrier}
	
	To clarify why the classical \(\mathrm{AT}_1\) model cannot independently prescribe material strength, we examine the onset of crack nucleation in a homogeneous intact body. The strong form of the phase-field evolution follows from the first-order optimality conditions of \cref{eq:local_min_continuous}. Taking the first variation of $\mathcal{E}_{\mathrm{AT}_1}^*$ with respect to $d$ gives the Karush--Kuhn--Tucker (KKT) conditions
	\begin{equation}
		\label{eq:kkt_classical}
		\dot{d} \geq 0,
		\qquad
		-g'(d)\psi_e^+
		-
		\frac{3G_c}{8l_0}
		+
		\frac{3G_c l_0}{4}\nabla^2 d
		\leq 0,
	\end{equation}
	together with the complementarity condition
	\begin{equation}
		\label{eq:complementarity_classical}
		\dot{d}
		\left[
		g'(d)\psi_e^+
		+
		\frac{3G_c}{8l_0}
		-
		\frac{3G_c l_0}{4}\nabla^2 d
		\right]
		=
		0.
	\end{equation}
	When damage grows, the constraint becomes active and the bracketed term vanishes:
	\begin{equation}
		\label{eq:active_damage_growth}
		\dot{d} > 0,
		\qquad
		g'(d)\psi_e^+
		+
		\frac{3G_c}{8l_0}
		-
		\frac{3G_c l_0}{4}\nabla^2 d
		=
		0.
	\end{equation}
	For an initially intact homogeneous body under monotonically increasing loading, we have $d=0$, $\nabla d=\mathbf{0}$, and thus $\nabla^2 d=0$ before nucleation. At the onset of damage, the active condition in \cref{eq:active_damage_growth} gives
	\begin{equation}
		g'(0)\psi_e^+
		+
		\frac{3G_c}{8l_0}
		=
		0.
	\end{equation}
	With $g(d)=(1-d)^2$, one has $g'(d)=-2(1-d)$ and $g'(0)=-2$. Therefore, the critical active elastic energy density required for damage initiation is
	\begin{equation}
		\label{eq:inherent_barrier}
		\psi_c
		=
		\frac{3G_c}{16l_0}.
	\end{equation}
	
	\begin{remark}[The intrinsic energy barrier and strength coupling]
		\label{rem:inherent_barrier}
		\Cref{eq:inherent_barrier} shows that the \(\mathrm{AT}_1\) model contains an intrinsic energy barrier for damage initiation. This barrier is controlled by the critical energy release rate $G_c$ and the regularization length $l_0$.
		
		For a one-dimensional tensile bar with Young's modulus $E$, the corresponding critical stress is
		$\sigma_c = \sqrt{3 G_c E / (8 l_0)}$.
		Thus, the nucleation stress is tied to $G_c$ and $l_0$. As a result, the tensile strength $\sigma_t$ and the critical energy release rate $G_c$ cannot be prescribed independently unless $l_0$ is adjusted accordingly. If $l_0$ is selected mainly from mesh-resolution considerations, the predicted failure load may become length-scale dependent.
	\end{remark}
	
	This coupling between strength, toughness, and the regularization length limits the direct use of the classical \(\mathrm{AT}_1\) model for materials whose failure is governed by a prescribed tensile strength. The following section introduces a shifted driving force to align the phase-field nucleation condition with the Rankine criterion while retaining the \(\mathrm{AT}_1\) fracture energy.

	\section{Shifted energy barrier approach for tensile-dominated fracture}
	\label{sec:shifted_theory}
	
	For brittle materials such as glass, ceramics, and concrete, crack initiation is often governed by tensile stresses. 
	A simple macroscopic criterion for tensile-dominated failure is the Rankine criterion, also known as the maximum principal stress criterion. We introduce the criterion into phase-field fracture in this section.

	\subsection{The Rankine criterion and critical nucleation energy}
	\label{subsec:rankine_energy}

	Let $\bm{\sigma}_0(\bm{\varepsilon})$ denote the effective undamaged elastic Cauchy stress tensor. 
	The principal stresses of $\bm{\sigma}_0$ are ordered as $\sigma_1\geq\sigma_2\geq\sigma_3$. 
	In the present formulation, crack nucleation is linked to the effective stress state and is assumed to occur when the maximum effective principal tensile stress reaches the uniaxial tensile strength $\sigma_t$:
	\begin{equation}
		\label{eq:rankine_initiation}
		\mathcal{F}(\bm{\sigma}_0)
		=
		\sigma_1
		-
		\sigma_t
		=
		0 .
	\end{equation}
	
	\begin{remark}
		The Rankine condition is evaluated using the effective stress $\bm{\sigma}_0$, rather than the degraded nominal stress. 
		This separates the local nucleation criterion from the subsequent stiffness degradation governed by the phase-field variable. 
		After damage initiation, the macroscopic stress response is still determined by the degraded stress, namely $\bm{\sigma}=g(d)\bm{\sigma}_0^+ + \bm{\sigma}_0^-$ under the adopted tension--compression split.
	\end{remark}
	
	To incorporate the stress-based Rankine criterion into the energy-driven \(\mathrm{AT}_1\) framework, the tensile strength $\sigma_t$ is mapped to an equivalent active elastic energy threshold, denoted by $\psi_{cs}$. 
	
	In general, however, a strength surface cannot be expressed as a \emph{fixed} energy or strain threshold; equivalent energy or strain thresholds exist only in special cases such as uniaxial tension \citep{khayaz2025comparison}. The threshold \(\psi_{cs}\) is therefore not a fixed energy level but is made state-dependent: since the active elastic energy of a general multiaxial state depends on the full strain tensor, \(\psi_{cs}\) is obtained by scaling the current effective stress state to the Rankine surface.
	
	For a given undamaged elastic state with $\sigma_1>0$, consider a virtual proportional scaling of the effective stress tensor until the condition $\sigma_1=\sigma_t$ is reached. 
	The corresponding scaling factor is
	\begin{equation}
		\eta
		=
		\frac{\sigma_t}{\sigma_1}.
	\end{equation}
	Since the elastic strain energy is quadratic in stress under linear elasticity, the active elastic energy scales with $\eta^2$. 
	The Rankine-equivalent critical active energy density is defined as
	\begin{equation}
		\label{eq:psi_cs}
		\psi_{cs}(\bm{\varepsilon})
		=
		\begin{cases}
			\eta^2
			\psi_e^+(\bm{\varepsilon}),
			& \text{if } \sigma_1>0, \\[6pt]
			\infty,
			& \text{if } \sigma_1\leq0 .
		\end{cases}
	\end{equation}
	The infinite threshold assigned for $\sigma_1\leq0$ precludes tensile crack
	nucleation under compressive effective stress states; its finite numerical
	surrogate is discussed in \cref{rem:compressive_threshold}.
	
	\begin{remark}[On non-proportional loading paths]
		The proportional scaling argument is used only to define an instantaneous
		energy threshold associated with the current effective stress state.
		It does not impose a proportional loading history on the actual deformation
		process: while a material point remains intact, $\psi_{cs}$ is re-evaluated
		from the current effective stress, so arbitrary, including non-proportional,
		loading paths are admissible in the elastic regime.
		The equivalence between the threshold $\psi_{cs}$ and the Rankine criterion
		at damage initiation is established in \cref{subsec:shifted_driving_force},
		and the treatment of the threshold after initiation is specified in
		\cref{sec:thermodynamic_admissibility}.
	\end{remark}

	\subsection{The shifted phase-field driving force}
	\label{subsec:shifted_driving_force}
	
	As shown in \cref{subsec:inherent_barrier}, the classical \(\mathrm{AT}_1\) model contains the intrinsic energy barrier $\psi_c=3G_c/(16l_0)$. 
	Damage starts when the active elastic energy $\psi_e^+$ reaches this barrier. 
	Since $\psi_c$ is controlled by $G_c$ and $l_0$, it does not generally coincide with the energy level associated with a prescribed tensile strength. 
	To link crack nucleation to the Rankine criterion, the damage driving term is shifted so that the onset condition becomes $\psi_e^+=\psi_{cs}$.
	
	We define the shifted active energy density as
	\begin{equation}
		\label{eq:shifted_driving_force}
		\widehat{\psi}_e^+
		(\bm{\varepsilon},\psi_{cs})
		=
		\psi_e^+(\bm{\varepsilon})
		-
		\psi_{cs}(\bm{\varepsilon})
		+
		\psi_c .
	\end{equation}
	This shift modifies the damage driving force, while the \(\mathrm{AT}_1\) crack density is kept unchanged.
	
	Replacing $\psi_e^+$ by $\widehat{\psi}_e^+$ in the active damage equation gives
	\begin{equation}
		\label{eq:shifted_evolution}
		g'(d)\widehat{\psi}_e^+
		+
		\frac{3G_c}{8l_0}
		-
		\frac{3G_c l_0}{4}\nabla^2 d
		=
		0 .
	\end{equation}
	At crack nucleation in an initially intact homogeneous state, $d=0$, $\nabla d=\mathbf 0$, and $\nabla^2d=0$. 
	Using $g'(0)=-2$, \cref{eq:shifted_evolution} gives
	\begin{equation}
		-2\widehat{\psi}_e^+
		+
		\frac{3G_c}{8l_0}
		=
		0 .
	\end{equation}
	Substituting \cref{eq:shifted_driving_force} and $\psi_c=3G_c/(16l_0)$ yields
	\begin{equation}
		\label{eq:shifted_onset_balance}
		-2
		\left(
		\psi_e^+
		-
		\psi_{cs}
		+
		\psi_c
		\right)
		+
		2\psi_c
		=
		0,
	\end{equation}
	and therefore
	\begin{equation}
		\label{eq:onset_condition}
		\psi_e^+(\bm{\varepsilon})
		=
		\psi_{cs}(\bm{\varepsilon}) .
	\end{equation}
	Thus, the shifted driving force changes the \(\mathrm{AT}_1\) nucleation
	condition from the intrinsic barrier $\psi_c$ to the Rankine-equivalent
	threshold $\psi_{cs}$.
	For $\sigma_1>0$, substituting \cref{eq:psi_cs} into \cref{eq:onset_condition}
	gives
	\[
	\left(1-\eta^2\right)\psi_e^+
	=
	\left[1-\left(\frac{\sigma_t}{\sigma_1}\right)^2\right]\psi_e^+
	=
	0 ,
	\]
	which, for $\psi_e^+>0$ with $\sigma_t>0$ and $\sigma_1>0$, reduces to
	$\sigma_1=\sigma_t$.
	The onset condition \cref{eq:onset_condition} is therefore equivalent to the
	Rankine criterion \cref{eq:rankine_initiation} evaluated at the current
	effective stress state.
	
	\begin{remark}[Compressive stress states]
		\label{rem:compressive_threshold}
		For $\sigma_1\leq0$, the threshold is set to $\psi_{cs}=\infty$ in
		\cref{eq:psi_cs} to preclude tensile nucleation under compression. In the
		numerical implementation, $\psi_{cs}$ is instead assigned the finite value
		$\beta\,\psi_c$. The driving force $\psi_e^+-\psi_{cs}+\psi_c$ remains
		negative whenever $\beta\,\psi_c$ exceeds the active energy $\psi_e^+$
		reached in the compressed region, which fixes the scale of $\beta$ relative
		to the ratio $\psi_e^+/\psi_c$. Any $\beta$ above this bound gives the same
		response, and we use $\beta=10^{10}$.
	\end{remark}
	
	After damage initiation, the current threshold $\psi_{cs}$ is replaced by a frozen history threshold, which will be introduced in \cref{sec:thermodynamic_admissibility}. 
	This avoids changes of the activated damage resistance during unloading or non-proportional changes in the stress state.
	
	\subsection{Restricted variational principle, stress formulation and KKT conditions}
	\label{subsec:variational_kkt}
	
	The classical variational phase-field formulation of \cref{sec:classical} derives from a single energy functional, so that the stress response and the damage driving force are governed by one common potential; the formulation is in this sense variationally consistent, with a self-adjoint structure \citep{Francfort1998,Bourdin2000,Pham2011a,wick2020multiphysics}.
	In the present shifted formulation this is no longer the case: because the Rankine-equivalent threshold \(\psi_{cs}\) depends on the strain through the effective stress, the equilibrium and damage equations cannot be obtained from a common potential.
	We therefore derive them from a restricted variational principle \citep{rosen1954solution,rosen1954use}, in which a state-dependent quantity is held fixed during the variation and its dependence on the primary fields is restored only afterwards.
	As analyzed by \citet{finlayson1967search}, the operator associated with such formulations is in general not self-adjoint, so they do not constitute genuine minimization principles. In this work, the restricted variational principle is used only as a systematic derivation device, the physical admissibility of the resulting evolution being established separately in \cref{sec:thermodynamic_admissibility}.
	
	In the present formulation, the state-dependent quantity treated in this restricted sense is the Rankine-equivalent threshold \(\psi_{cs}=\psi_{cs}(\bm{\varepsilon}(\mathbf u))\).
	At a given load step, \(\psi_{cs}\) is evaluated from the effective stress state used to define the damage-initiation threshold and is then kept fixed during the variations with respect to the primary fields.
	We introduce the restricted shifted functional
	\begin{equation}
		\label{eq:shifted_functional}
		\begin{aligned}
			\mathcal{E}_{\mathrm{shift}}(\mathbf u,d;\psi_{cs})
			=&
			\int_{\Omega}
			\left[
			g(d)
			\left(
			\psi_e^+(\bm{\varepsilon}(\mathbf u))
			-
			\psi_{cs}
			+
			\psi_c
			\right)
			+
			\psi_e^-(\bm{\varepsilon}(\mathbf u))
			\right]d\Omega  \\
			&+
			\frac{3G_c}{8l_0}
			\int_{\Omega}
			\left(
			d+l_0^2|\nabla d|^2
			\right)d\Omega .
		\end{aligned}
	\end{equation}
	The restricted shifted functional is used only to derive the mechanical equilibrium equation and the shifted phase-field conditions under this prescribed-threshold restriction; it is not the Helmholtz free energy entering the thermodynamic dissipation inequality.
	The passive elastic energy \(\psi_e^-\) is retained in the displacement variation to preserve the standard compressive elastic response, but does not contribute to the phase-field driving force because it is not degraded by \(g(d)\).
	
	Consider first the variation with respect to the displacement field.
	Since \(\psi_{cs}\) depends on the strain through the effective stress, a full variation of \(\mathcal{E}_{\mathrm{shift}}\) would yield the stress response
	\begin{equation}
		\label{eq:full_variation_stress}
		\bm{\sigma}_{\mathrm{full}}
		=
		g(d)
		\left(
		\bm{\sigma}_0^+
		-
		\frac{\partial \psi_{cs}}{\partial\bm{\varepsilon}}
		\right)
		+
		\bm{\sigma}_0^- ,
	\end{equation}
	in which the additional term \(-g(d)\,\partial\psi_{cs}/\partial\bm{\varepsilon}\) would modify the standard degraded elastic response and is not part of the constitutive assumptions of the present model.
	In the restricted variation, \(\psi_{cs}\) (together with the constant \(\psi_c\)) is held fixed, so its derivative with respect to strain vanishes and the stress reduces to
	\begin{equation}
		\label{eq:shifted_stress}
		\bm{\sigma}(\bm{\varepsilon},d)
		=
		\left.
		\frac{\partial}{\partial\bm{\varepsilon}}
		\left[
		g(d)
		\left(
		\psi_e^+
		-
		\psi_{cs}
		+
		\psi_c
		\right)
		+
		\psi_e^-
		\right]
		\right|_{\psi_{cs}\ \text{fixed}}
		=
		g(d)\bm{\sigma}_0^+
		+
		\bm{\sigma}_0^- .
	\end{equation}
	Thus the restricted variation lets the shifted threshold modify the phase-field driving force while leaving the degraded elastic stress unchanged.
	This is precisely the mechanism by which the restricted variational principle provides a route to a hybrid formulation \citep{Ambati2015review}: the prescribed threshold drives the damage evolution without altering the mechanical equilibrium.
	The use of a non-energetic, strength-based contribution to the damage driving force is shared by recent phase-field formulations that introduce material strength as an independent input \citep{Kumar2020,Larsen2024,LopezPamies2025}; here it is realized through the restricted variation of the shifted threshold while the standard \(\mathrm{AT}_1\) stress is retained.
	
	Taking the restricted variation with respect to the phase-field variable \(d\) gives
	\begin{equation}
		\label{eq:restricted_damage_derivative}
		\frac{\delta \mathcal{E}_{\mathrm{shift}}}{\delta d}
		=
		g'(d)
		\left(
		\psi_e^+
		-
		\psi_{cs}
		+
		\psi_c
		\right)
		+
		\frac{3G_c}{8l_0}
		-
		\frac{3G_c l_0}{4}
		\nabla^2 d .
	\end{equation}
	The irreversible damage evolution is then governed by the KKT (loading--unloading) conditions associated with the irreversibility constraint \(\dot d\ge0\) and the driving force in \cref{eq:restricted_damage_derivative}:
	\begin{subnumcases}{}
		\dot{d} \geq 0,
		\label{eq:kkt_primal} \\[2pt]
		-g'(d)
		\left(
		\psi_e^+
		-
		\psi_{cs}
		+
		\psi_c
		\right)
		-
		\frac{3G_c}{8l_0}
		+
		\frac{3G_c l_0}{4}
		\nabla^2 d
		\leq 0,
		\label{eq:kkt_dual} \\[2pt]
		\dot{d}
		\left[
		g'(d)
		\left(
		\psi_e^+
		-
		\psi_{cs}
		+
		\psi_c
		\right)
		+
		\frac{3G_c}{8l_0}
		-
		\frac{3G_c l_0}{4}
		\nabla^2 d
		\right]
		=
		0 .
		\label{eq:kkt_comp}
	\end{subnumcases}
	These conditions are stated as a complementarity (variational-inequality) problem rather than as the optimality conditions of a global minimization. They coincide with the optimality conditions of the bound-constrained phase-field subproblem solved at held-fixed threshold (\cref{sec:numerical_implementation}), and the driving force entering them is the one obtained thermodynamically from the dissipative microforce balance in \cref{sec:thermodynamic_admissibility}, where the associated damage dissipation is shown to be non-negative.
	
	For an initially intact homogeneous state, \(d=0\) and
	\(\nabla^2 d=0\), so the active complementarity condition
	\cref{eq:kkt_comp} reduces to \cref{eq:shifted_onset_balance},
	and the nucleation condition \(\psi_e^+=\psi_{cs}\) of
	\cref{eq:onset_condition} is recovered. The restricted variational
	principle therefore reproduces the shifted nucleation condition
	while preserving the standard degraded elastic stress response in
	\cref{eq:shifted_stress}.
	
	\section{Thermodynamics}
	\label{sec:thermodynamic_admissibility}
	
	The restricted variational formulation in \cref{subsec:variational_kkt} gives the shifted phase-field equation while preserving the degraded elastic stress response in \cref{eq:shifted_stress}. 
	In the restricted variation, the Rankine-equivalent threshold \(\psi_{cs}\) is treated as a fixed local parameter during variation. 
	Therefore, the threshold enters the phase-field equation but not the displacement equilibrium equation.
	
	The thermodynamic interpretation adopted here keeps the recoverable Helmholtz free energy equal to the standard \(\mathrm{AT}_1\) density associated with \cref{eq:at1_total_energy}. 
	The shifted functional in \cref{eq:shifted_functional} is not regarded as an independent recoverable free energy. 
	Instead, the Rankine-based shift is interpreted as a dissipative microforce contribution conjugate to the irreversible evolution of the phase-field variable. 
	The formulation follows microforce-based thermodynamics for gradient-type internal variables \citep{gurtin1996generalized,duda2015phase}.
	
	\subsection{Microforce-based thermodynamics}
	\label{subsec:microforce_thermo_setting}
	
	The recoverable Helmholtz free-energy density is taken as the standard
	\(\mathrm{AT}_1\) density, i.e.\ the integrand of
	\cref{eq:at1_total_energy},
	\begin{equation}
		\label{eq:Psi0_definition}
		\Psi_0(\bm{\varepsilon},d,\nabla d)
		=
		g(d)\,\psi_e^+(\bm{\varepsilon})
		+
		\psi_e^-(\bm{\varepsilon})
		+
		\frac{3G_c}{8l_0}
		\left(
		d
		+
		l_0^2\,|\nabla d|^2
		\right),
	\end{equation}
	so that
	\(\mathcal{E}_{\mathrm{AT}_1}^*=\int_{\Omega}\Psi_0\,d\Omega\).
	Since the Rankine-equivalent threshold \(\psi_{cs}\) is not included in
	\(\Psi_0\), the recoverable stress remains the degraded elastic stress
	already given in \cref{eq:shifted_stress}.
	
	Because \(\Psi_0\) depends on \(d\) and \(\nabla d\), the phase-field variable is treated as a gradient-type internal variable. 
	A scalar microforce \(\pi\), conjugate to \(\dot d\), and a vector microstress \(\bm{\xi}\), conjugate to \(\nabla\dot d\), are introduced. 
	The local isothermal free-energy imbalance is written as
	\begin{equation}
		\label{eq:microforce_dissipation}
		\mathcal D
		=
		\bm{\sigma}:\dot{\bm{\varepsilon}}
		+
		\pi \dot d
		+
		\bm{\xi}\cdot\nabla \dot d
		-
		\dot{\Psi}_0
		\ge 0 .
	\end{equation}
	Application of the Coleman--Noll procedure \citep{coleman1963thermodynamics} to the recoverable rates gives
	\begin{equation}
		\label{eq:recoverable_stress_microstress}
		\bm{\sigma}
		=
		\frac{\partial\Psi_0}{\partial\bm{\varepsilon}},
		\qquad
		\bm{\xi}
		=
		\frac{\partial\Psi_0}{\partial\nabla d}.
	\end{equation}
	The remaining contribution to the dissipation is associated with the scalar microforce conjugate to \(\dot d\). 
	We define
	\begin{equation}
		\label{eq:pi_dis_definition}
		\pi^{\mathrm{dis}}
		:=
		\pi-\frac{\partial\Psi_0}{\partial d},
		\qquad
		\mathcal D_d
		=
		\pi^{\mathrm{dis}}\dot d
		\ge0 .
	\end{equation}
	The coefficient of \(\dot d\) in \cref{eq:pi_dis_definition} is not set to zero. Crack irreversibility restricts the admissible damage rates to \(\dot d\ge0\).
	Therefore, the remaining scalar microforce contribution is dissipative.

	In the absence of external microforces associated with \(d\), the local microforce balance is
	\begin{equation}
		\label{eq:microforce_balance}
		\pi-\nabla\cdot\bm{\xi}=0 .
	\end{equation}
	Combining \cref{eq:recoverable_stress_microstress,eq:pi_dis_definition,eq:microforce_balance} gives
	\begin{equation}
		\label{eq:microforce_balance_expanded}
		\frac{\partial\Psi_0}{\partial d}
		-
		\nabla\cdot
		\left(
		\frac{\partial\Psi_0}{\partial\nabla d}
		\right)
		+
		\pi^{\mathrm{dis}}
		=
		0 .
	\end{equation}
	For the standard \(\mathrm{AT}_1\) density, the energetic terms in \cref{eq:microforce_balance_expanded} recover the classical phase-field balance derived in \cref{subsec:inherent_barrier}. 
	Consequently, a Rankine-based modification of the phase-field evolution must enter through the dissipative scalar microforce, rather than through the recoverable free energy.
	
	\subsection{Rankine-based dissipative microforce}
	\label{subsec:rankine_based_dissipative_microforce}
	
	Before damage initiation, the threshold controlling the shifted response is the current Rankine-equivalent value \(\psi_{cs}\), because the intended onset condition is \(\psi_e^+=\psi_{cs}\). 
	After damage initiation, continuous updating of \(\psi_{cs}\) from the current stress state would make the activated resistance vary during unloading or non-proportional loading. 
	The threshold entering the dissipative microforce is therefore frozen after initiation:
	\begin{equation}
		\label{eq:Hcs_freezing_thermo}
		\mathcal H_{cs}(\mathbf{x},t)
		=
		\begin{cases}
			\psi_{cs}(\bm{\varepsilon}(\mathbf{x},t)),
			& d(\mathbf{x},t)=0, \\[4pt]
			\psi_{cs}(\bm{\varepsilon}(\mathbf{x},t_{\mathrm{ini}})),
			& d(\mathbf{x},t)>0 ,
		\end{cases}
	\end{equation}
	where \(t_{\mathrm{ini}}(\mathbf{x})\) is the local damage-initiation time, defined as the first instant at which \(d(\mathbf{x},t)>0\). Accordingly, \(\mathcal H_{cs}\) follows the current Rankine-equivalent threshold while the material point is intact and is held fixed at its initiation value once damage has started, so that \(\dot{\mathcal H}_{cs}=0\) for \(d>0\).
	In this sense \(\mathcal H_{cs}\) is the central quantity of the shifted formulation: before initiation it fixes the strength-controlled nucleation condition \(\psi_e^+=\mathcal H_{cs}\), and once frozen it sets the resistance to subsequent damage growth. Because it is frozen rather than load-following, the Rankine shift enters as a well-defined dissipative resistance, whose thermodynamic admissibility is established in \cref{subsec:dissipation_admissibility_reduced}.
	
	The dissipative scalar microforce associated with the Rankine-based shift is prescribed as
	\begin{equation}
		\label{eq:pi_dis_rankine}
		\pi^{\mathrm{dis}}
		=
		-
		g'(d)
		\left(
		\mathcal H_{cs}
		-
		\psi_c
		\right),
	\end{equation}
	where the intrinsic \(\mathrm{AT}_1\) barrier \(\psi_c\) is given in \cref{eq:inherent_barrier}. 
	The factor \(-g'(d)\) makes the additional resistance enter the phase-field balance with the same degradation-weighted structure as the active elastic driving contribution.
	
	Substitution of \cref{eq:pi_dis_rankine} into \cref{eq:microforce_balance_expanded} gives the phase-field balance associated with active damage growth:
	\begin{equation}
		\label{eq:shifted_balance_thermo}
		g'(d)
		\left(
		\psi_e^+
		-
		\mathcal H_{cs}
		+
		\psi_c
		\right)
		+
		\frac{3G_c}{8l_0}
		-
		\frac{3G_c l_0}{4}
		\nabla^2 d
		=
		0 .
	\end{equation}
	The balance equation in \cref{eq:shifted_balance_thermo} has the same residual form as the equation obtained from the restricted variation in \cref{eq:restricted_damage_derivative}, with \(\psi_{cs}\) replaced by the frozen threshold \(\mathcal H_{cs}\) after initiation. 
	The complete irreversible evolution remains governed by the constrained phase-field problem and the KKT conditions stated in \cref{subsec:variational_kkt}. 
	The role of the thermodynamic argument is to identify the shifted term as a dissipative microforce contribution and to examine the sign of the associated damage dissipation.
	
	Before damage initiation, \(\mathcal H_{cs}=\psi_{cs}\). 
	For an initially intact state, \cref{eq:shifted_balance_thermo} reduces to the same nucleation condition \(\psi_e^+=\psi_{cs}\). 
	As shown in \cref{subsec:shifted_driving_force,subsec:variational_kkt}, this condition is equivalent to the Rankine criterion on the tensile branch, where the Rankine-equivalent threshold is finite.
	The dissipative-microforce formulation therefore preserves the Rankine-controlled nucleation condition derived from the restricted variational formulation.
	
	\begin{remark}
		Since the energy threshold is frozen as a scalar history variable upon damage initiation ($d>0$), the current formulation is limited to monotonic loading paths without significant rotation of principal stress directions. Arbitrary non-proportional loading paths are fully supported only in the preceding elastic regime ($d=0$).
	\end{remark}
	
	\subsection{Dissipation admissibility}
	\label{subsec:dissipation_admissibility_reduced}
	
	Using \cref{eq:pi_dis_rankine}, the local damage dissipation becomes
	\begin{equation}
		\label{eq:damage_dissipation_rankine}
		\mathcal D_d
		=
		-
		g'(d)
		\left(
		\mathcal H_{cs}
		-
		\psi_c
		\right)
		\dot d .
	\end{equation}
	Because irreversibility requires \(\dot d\ge0\), the case \(\dot d=0\) gives zero dissipation. 
	For all admissible damage growth processes, non-negative dissipation is guaranteed if
	\begin{equation}
		\label{eq:Hcs_admissibility}
		\mathcal H_{cs}
		\ge
		\psi_c ,
	\end{equation}
	where \(g(d)=(1-d)^2\) has been used. 
	Since \(\mathcal H_{cs}\) is the frozen value of \(\psi_{cs}\) at initiation, \cref{eq:Hcs_admissibility} requires \(\psi_{cs}(t_{\mathrm{ini}})\ge\psi_c\). 
	For uniaxial tension, the Rankine-equivalent threshold reduces to the constant
	value \(\psi_{cs}=\sigma_t^2/(2E)\), as shown in the one-dimensional analysis
	of \cref{sec:1d_analysis} (see \cref{eq:1d_threshold}). 
	Together with \(\psi_c=3G_c/(16l_0)\), the admissibility condition gives
	\begin{equation}
		\label{eq:l0_condition_thermo}
		l_0
		\ge
		\frac{3EG_c}{8\sigma_t^2}.
	\end{equation}
	The bound in \cref{eq:l0_condition_thermo} is not a mesh-size requirement, but the condition under which the Rankine-based shift can be interpreted as a non-negative dissipative resistance. 
	If the Rankine-equivalent threshold is lower than the intrinsic \(\mathrm{AT}_1\) barrier, the shift reduces the resistance to damage growth, and the sign of the dissipation in \cref{eq:damage_dissipation_rankine} is not guaranteed.
	
	\section{Softening problem of a one-dimensional bar}
	\label{sec:1d_analysis}
	
	The shifted formulation of \cref{sec:shifted_theory,sec:thermodynamic_admissibility} is constructed at the level of the multiaxial effective stress state. When it is specialized to a bar under uniaxial traction, the Rankine-equivalent threshold reduces to the constant energy level \(\sigma_t^2/(2E)\) and the freezing rule of \cref{eq:Hcs_freezing_thermo} has no effect. The one-dimensional problem therefore does not exercise the state dependence of the threshold, but it admits closed-form solutions, which quantify the effect of the barrier shift on the homogeneous softening response and on the localized damage profile and serve as reference solutions for the verification of the numerical implementation in \cref{subsec:uniaxial_tension}. Both solutions are derived following the classical constructions for gradient-damage models \citep{Pham2011a,Pham2011b,pham2013onset}, and the localized profile is compared with that of the classical \(\mathrm{AT}_1\) model.
	Throughout this section, \((\cdot)'\) denotes differentiation with respect to the axial coordinate \(x\).
	
	\subsection{One-dimensional reduction and the shift parameter}
	\label{subsec:1d_setting}
	
	Consider a bar under uniaxial tension with axial strain \(\varepsilon\ge0\) and Young's modulus \(E\).
	For \(n=1\), the deviatoric strain vanishes, \(\bm{\varepsilon}_{\mathrm{dev}}=\mathbf 0\), and Amor's decomposition in \cref{eq:amor_split} reduces to
	\begin{equation}
		\label{eq:1d_split}
		\psi_e^+(\varepsilon)
		=
		\frac{E}{2}\langle\varepsilon\rangle_+^2,
		\qquad
		\psi_e^-(\varepsilon)
		=
		\frac{E}{2}\langle\varepsilon\rangle_-^2 ,
	\end{equation}
	where the one-dimensional modulus \(E\) plays the role of \(\kappa\).
	Under tension, the effective stress is \(\sigma_0=E\varepsilon=\sigma_1>0\), and the Rankine-equivalent threshold in \cref{eq:psi_cs} becomes
	\begin{equation}
		\label{eq:1d_threshold}
		\psi_{cs}
		=
		\left(\frac{\sigma_t}{E\varepsilon}\right)^2
		\frac{E}{2}\varepsilon^2
		=
		\frac{\sigma_t^2}{2E},
	\end{equation}
	which is constant and independent of the strain state.
	The freezing rule in \cref{eq:Hcs_freezing_thermo} therefore has no effect on the damage evolution. The threshold equals \(\sigma_t^2/(2E)\) at every intact point under tension and is frozen at the same value wherever damage has initiated, so that
	\begin{equation}
		\label{eq:1d_Hcs}
		\mathcal H_{cs}
		=
		\frac{\sigma_t^2}{2E}
		\qquad
		\text{wherever } \sigma_0>0 \text{ or } d>0 .
	\end{equation}
	
	The shift with respect to the intrinsic \(\mathrm{AT}_1\) barrier \(\psi_c=3G_c/(16l_0)\) in \cref{eq:inherent_barrier} is measured by the dimensionless parameter
	\begin{equation}
		\label{eq:gamma_definition}
		\gamma
		:=
		\frac{\mathcal H_{cs}-\psi_c}{\psi_c}
		=
		\frac{8\sigma_t^2 l_0}{3EG_c}
		-
		1
		=
		\left(\frac{\sigma_t}{\sigma_c}\right)^2
		-
		1 ,
	\end{equation}
	where \(\sigma_c=\sqrt{3G_cE/(8l_0)}\) is the nucleation stress of the classical \(\mathrm{AT}_1\) model (\cref{rem:inherent_barrier}).
	The admissibility condition \(\mathcal H_{cs}\ge\psi_c\) in \cref{eq:Hcs_admissibility,eq:l0_condition_thermo} is equivalent to \(\gamma\ge0\).
	Thus, \(\gamma\) quantifies the amount by which the prescribed strength exceeds the intrinsic \(\mathrm{AT}_1\) nucleation stress at the given regularization length, and the limit \(\gamma\to0^+\) corresponds to the admissibility boundary.
	
	Using \(3G_c/(8l_0)=2\psi_c\), \(3G_cl_0/4=4\psi_cl_0^2\), and \(\mathcal H_{cs}-\psi_c=\gamma\,\psi_c\), the shifted phase-field balance \cref{eq:shifted_balance_thermo} for active damage growth (\(\dot d>0\)) takes the one-dimensional form
	\begin{equation}
		\label{eq:1d_balance}
		-2(1-d)
		\left(
		\frac{E}{2}\varepsilon^2
		-
		\gamma\,\psi_c
		\right)
		+
		2\psi_c
		-
		4\psi_c l_0^2\, d''
		=
		0 .
	\end{equation}
	
	\begin{remark}[Constant threshold in one dimension]
		\label{rem:1d_constant_threshold}
		In the uniaxial setting the shifted formulation thus coincides with an \(\mathrm{AT}_1\) model whose damage driving force is translated by a constant. The state dependence of \(\psi_{cs}\), and hence the role of the frozen threshold \(\mathcal H_{cs}\), appears only under multiaxial stress states, where no fixed energy threshold equivalent to the Rankine criterion exists \citep{khayaz2025comparison}. These ingredients are assessed in \cref{subsec:hole_compression,subsec:multiaxial,subsec:cruciform_2d}.
	\end{remark}
	
	\subsection{Homogeneous solution}
	\label{subsec:1d_homogeneous}
	
	Consider a homogeneous state \((\varepsilon,d)\) with \(d''=0\) under monotonically increasing strain.
	The intact state \(d=0\) is admissible as long as the one-dimensional form of the KKT inequality \cref{eq:kkt_dual} holds,
	\begin{equation}
		\label{eq:1d_elastic_domain}
		2\left(\frac{E}{2}\varepsilon^2-\mathcal H_{cs}\right)\le0
		\quad\Longleftrightarrow\quad
		\sigma_0=E\varepsilon\le\sigma_t .
	\end{equation}
	Damage therefore initiates at the Rankine limit \(\varepsilon_t=\sigma_t/E\), independently of \(l_0\), consistent with \cref{subsec:shifted_driving_force}.
	
	For \(\varepsilon\ge\varepsilon_t\), the consistency condition \cref{eq:1d_balance} gives the damaging branch
	\begin{equation}
		\label{eq:1d_homogeneous_d}
		1-d
		=
		\frac{\psi_c}{\dfrac{E}{2}\varepsilon^2-\gamma\,\psi_c} .
	\end{equation}
	The denominator is positive on this branch, since \(\tfrac{E}{2}\varepsilon^2\ge\mathcal H_{cs}=(1+\gamma)\psi_c>\gamma\,\psi_c\).
	At \(\varepsilon=\varepsilon_t\), the denominator equals \(\mathcal H_{cs}-\gamma\psi_c=\psi_c\), so that \(d=0\) and the elastic and damaging branches connect continuously.
	Moreover, \(d\) in \cref{eq:1d_homogeneous_d} increases monotonically with \(\varepsilon\) and tends to \(1\) as \(\varepsilon\to\infty\), so the irreversibility constraint \(\dot d\ge0\) is satisfied along monotone loading.
	
	Introducing the dimensionless effective stress \(\tilde{\sigma}=\sigma_0/\sigma_t=\varepsilon/\varepsilon_t\ge1\) and using \(\tfrac{E}{2}\varepsilon^2=(1+\gamma)\tilde{\sigma}^2\psi_c\), the damaging branch and the nominal (degraded) stress \(\sigma=(1-d)^2E\varepsilon\) become
	\begin{equation}
		\label{eq:1d_homogeneous_dimensionless}
		1-d
		=
		\frac{1}{(1+\gamma)\tilde{\sigma}^2-\gamma},
		\qquad
		\frac{\sigma}{\sigma_t}
		=
		\frac{\tilde{\sigma}}{\left[(1+\gamma)\tilde{\sigma}^2-\gamma\right]^{2}} ,
		\qquad
		\tilde{\sigma}\ge1 .
	\end{equation}
	For \(\gamma=0\), where \(\sigma_t=\sigma_c\), \cref{eq:1d_homogeneous_dimensionless} reduces to \(1-d=\tilde{\sigma}^{-2}\) and \(\sigma/\sigma_c=\tilde{\sigma}^{-3}\), which is the homogeneous softening response of the classical \(\mathrm{AT}_1\) model \citep{Pham2011a}.
	
	The peak nominal stress is attained at initiation.
	Writing \(\psi=\tfrac{E}{2}\varepsilon^2\), the nominal stress reads \(\sigma=\psi_c^2E\varepsilon/(\psi-\gamma\psi_c)^2\), and differentiation with respect to \(\varepsilon\) gives
	\begin{equation}
		\label{eq:1d_softening_slope}
		\frac{d\sigma}{d\varepsilon}
		=
		-\,
		\frac{\psi_c^2\,E\,\left(3\psi+\gamma\psi_c\right)}{\left(\psi-\gamma\psi_c\right)^{3}}
		<0
		\qquad
		\text{for }\varepsilon\ge\varepsilon_t .
	\end{equation}
	The homogeneous response therefore softens immediately after initiation, without a hardening plateau, and the peak nominal stress equals the prescribed strength \(\sigma_t\).
	A stability analysis of the homogeneous states in the sense of \citet{Pham2011b} is beyond the scope of this work.
	Since the branch \cref{eq:1d_homogeneous_dimensionless} softens immediately, however, localization at the peak is expected for sufficiently long bars, as in the classical \(\mathrm{AT}_1\) model \citep{Pham2011b,Zolesi2024}.
	
	\subsection{Localized solution and optimal damage profile}
	\label{subsec:1d_localized}
	
	The localized solution corresponds to the ultimate damage profile of a stress-free crack and is constructed following the classical procedure for gradient-damage models \citep{pham2013onset,Tanne2018}.
	Let the localization band occupy \(|x|\le D\), with the crack center at \(x=0\).
	At complete failure the axial stress vanishes, \(\sigma=(1-d)^2E\varepsilon=0\), so that \(\varepsilon=0\) and hence \(\psi_e^+=0\) at every point of the band where \(d<1\).
	The profile is sought symmetric and monotone, with
	\begin{equation}
		\label{eq:1d_profile_bcs}
		d(0)=1,
		\qquad
		d(\pm D)=0,
		\qquad
		d'(\pm D)=0,
	\end{equation}
	where the last condition ensures a smooth (\(C^1\)) matching with the undamaged outer region, and the balance \cref{eq:1d_balance} holds with equality on the support of \(d\).
	By \cref{eq:1d_Hcs}, the frozen threshold takes the same value \(\sigma_t^2/(2E)\) at every point of the band, irrespective of the local initiation time (\cref{rem:1d_constant_threshold}).
	
	Setting \(\psi_e^+=0\) in \cref{eq:1d_balance} and introducing \(e=1-d\), so that \(d''=-e''\), one obtains for \(\gamma>0\) the linear ordinary differential equation
	\begin{equation}
		\label{eq:1d_ode}
		e''
		+
		\omega^2 e
		=
		-\frac{1}{2l_0^2},
		\qquad
		\omega
		=
		\frac{1}{l_0}\sqrt{\frac{\gamma}{2}} ,
	\end{equation}
	on \(x\in[0,D]\), with \(e(0)=0\), \(e(D)=1\), and \(e'(D)=0\).
	For the classical \(\mathrm{AT}_1\) model, the corresponding equation is \(e''=-1/(2l_0^2)\).
	The shift \(\mathcal H_{cs}-\psi_c=\gamma\psi_c\) introduces the additional restoring term \(\omega^2e\), which changes the profile from parabolic to trigonometric type.
	
	The general solution of \cref{eq:1d_ode} is
	\begin{equation}
		e(x)
		=
		A\cos\omega x
		+
		B\sin\omega x
		-
		\frac{1}{\gamma},
	\end{equation}
	where the particular solution follows from \(1/(2l_0^2\omega^2)=1/\gamma\).
	The condition \(e(0)=0\) gives \(A=1/\gamma\), and the condition \(e'(D)=0\) gives \(B=A\tan\theta\), with \(\theta:=\omega D\).
	Substituting these constants into the remaining condition \(e(D)=1\) yields \(A/\cos\theta-1/\gamma=1\), and therefore
	\begin{equation}
		\label{eq:1d_theta}
		\cos\theta
		=
		\frac{1}{1+\gamma},
		\qquad
		\theta\in\left(0,\frac{\pi}{2}\right).
	\end{equation}
	The half-width of the localization band and the damage profile then admit the closed forms
	\begin{equation}
		\label{eq:1d_bandwidth}
		D
		=
		\frac{\theta}{\omega}
		=
		\frac{\sqrt2\,l_0}{\sqrt\gamma}
		\arccos\frac{1}{1+\gamma} ,
	\end{equation}
	\begin{equation}
		\label{eq:1d_profile}
		d(x)
		=
		\frac{1+\gamma}{\gamma}
		\left[
		1-\cos\left(\theta-\omega|x|\right)
		\right],
		\qquad
		|x|\le D,
	\end{equation}
	with \(d(x)=0\) for \(|x|\ge D\).
	
	The profile \cref{eq:1d_profile} satisfies \(d(0)=1\) by \cref{eq:1d_theta} and decreases monotonically in \(|x|\), since \(d'(x)\propto-\sin(\theta-\omega|x|)\le0\) for \(\theta-\omega|x|\in[0,\theta]\).
	Hence \(0\le d\le1\) throughout.
	Outside the band, \(d=0\) and \(\varepsilon=0\) give \(\sigma_1=0\), so that \(\psi_{cs}=\infty\) by \cref{eq:psi_cs} and the KKT inequality \cref{eq:kkt_dual} is satisfied.
	At the boundary of the support, the interior limit of the curvature is \(d''(D^-)=(1+\gamma)/(2l_0^2)>0\), with a jump to zero across \(x=D\).
	The classical \(\mathrm{AT}_1\) profile exhibits the same finite curvature jump at the boundary of its support, with \(d''=1/(2l_0^2)\), and both profiles are compatible with the \(H^1\) regularity of the phase field.
	
	\subsection{Comparison with the classical \(\mathrm{AT}_1\) profile}
	\label{subsec:1d_comparison}
	
	The localized profile of the classical \(\mathrm{AT}_1\) model is parabolic \citep{Pham2011a,Tanne2018},
	\begin{equation}
		\label{eq:at1_profile}
		d_{\mathrm{AT}_1}(x)
		=
		\left(1-\frac{|x|}{2l_0}\right)^2,
		\qquad
		|x|\le 2l_0 ,
	\end{equation}
	with half-width \(2l_0\).
	The shifted profile \cref{eq:1d_profile} is of cosine type, and its half-width \cref{eq:1d_bandwidth} depends on the shift parameter \(\gamma\).
	Two properties relate the two profiles.
	
	First, the classical profile is recovered in the limit \(\gamma\to0^+\).
	Expanding \cref{eq:1d_theta} for small \(\gamma\) gives
	\(\theta=\sqrt{2\gamma}\,\left(1-\tfrac{5}{12}\gamma+O(\gamma^2)\right)\),
	and therefore
	\begin{equation}
		\label{eq:1d_bandwidth_expansion}
		D
		=
		2l_0
		\left(
		1-\frac{5}{12}\,\gamma+O(\gamma^2)
		\right)
		\rightarrow
		2l_0
		\qquad
		\text{as }\gamma\to0^+ .
	\end{equation}
	A Taylor expansion of \cref{eq:1d_profile} for small \(\gamma\), with \(\theta\) and \(\omega|x|\) both of order \(\sqrt\gamma\), gives \(d(x)\to(1-|x|/(2l_0))^2\), so the parabolic profile \cref{eq:at1_profile} is recovered pointwise at the admissibility boundary.
	The slope at the crack center behaves consistently: \(|d'(0^\pm)|=\sqrt{(\gamma+2)/2}\,/\,l_0\to1/l_0\), which matches the kink of the \(\mathrm{AT}_1\) profile at \(x=0\).
	
	Second, the localization band is strictly narrower than the \(\mathrm{AT}_1\) band for every \(\gamma>0\),
	\begin{equation}
		\label{eq:1d_band_bound}
		D<2l_0 .
	\end{equation}
	To verify \cref{eq:1d_band_bound}, note that \(D<2l_0\) is equivalent to \(\theta<\sqrt{2\gamma}\) with \(\gamma=(1-\cos\theta)/\cos\theta\), and hence to
	\begin{equation}
		q(\theta)
		:=
		2(1-\cos\theta)-\theta^2\cos\theta
		>0
		\qquad
		\text{on }\left(0,\frac{\pi}{2}\right).
	\end{equation}
	Since \(q(0)=0\) and \(q'(\theta)=2(\sin\theta-\theta\cos\theta)+\theta^2\sin\theta>0\) on this interval, as \(\tan\theta>\theta\), the bound follows.
	In addition, \(D/l_0\) decreases monotonically with \(\gamma\), with the asymptotic behavior \(D\approx\pi l_0/\sqrt{2\gamma}\) as \(\gamma\to\infty\).

	\begin{figure}[htbp]
		\centering
		\includegraphics[width=1.0\textwidth]{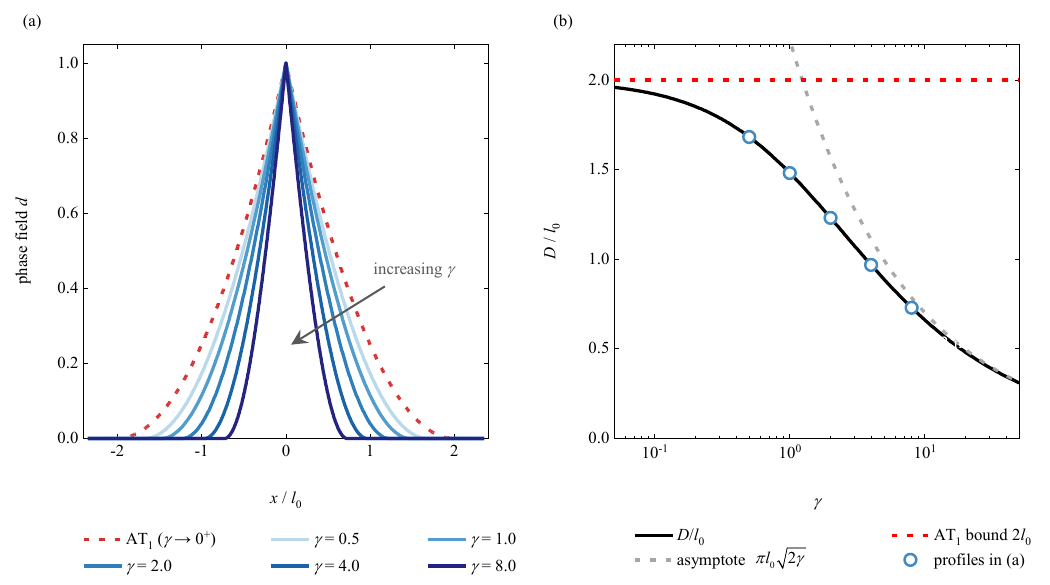}
		\caption{Ultimate localized damage profiles in one dimension.
			(a) Cosine-type profiles of the shifted formulation, \cref{eq:1d_profile}, for increasing values of the shift parameter \(\gamma\) defined in \cref{eq:gamma_definition}, together with the classical parabolic \(\mathrm{AT}_1\) profile \cref{eq:at1_profile} recovered in the limit \(\gamma\to0^+\).
			(b) Normalized half-width \(D/l_0\) of the localization band, \cref{eq:1d_bandwidth}, as a function of \(\gamma\): the half-width decreases monotonically, is bounded by the classical value \(2l_0\) (red dashed line), and approaches the asymptote \(\pi/\sqrt{2\gamma}\) of \cref{rem:1d_band_width_scaling} for large \(\gamma\) (gray dashed line).
			Open circles mark the values of \(\gamma\) corresponding to the profiles shown in (a).}
		\label{fig:1d_profiles}
	\end{figure}

	\begin{remark}[Band width at large shift]
		\label{rem:1d_band_width_scaling}
		For \(\gamma\gg1\), combining \(D\approx\pi l_0/\sqrt{2\gamma}\) with \cref{eq:gamma_definition} gives
		\[
		D
		\approx
		\frac{\pi\sqrt3}{4}
		\sqrt{l_0\, l_{\mathrm{ch}}},
		\qquad
		l_{\mathrm{ch}}
		:=
		\frac{EG_c}{\sigma_t^2},
		\]
		where \(l_{\mathrm{ch}}\) is the Irwin-type material length.
		At large regularization lengths, the band width therefore grows only like \(\sqrt{l_0}\), in contrast with the linear scaling \(2l_0\) of the classical \(\mathrm{AT}_1\) model: the shift partially compensates the widening of the regularized crack.
	\end{remark}

	\section{Numerical implementation}
	\label{sec:numerical_implementation}
	
	The coupled problem is discretized in space by the standard finite element 
	method and solved by an under-relaxed staggered scheme. 
	At time step \(t_{n+1}\), the iteration is initialized from the converged state at \(t_n\), namely 
	\(\mathbf u^{(0)}=\mathbf u_n\), \(d^{(0)}=d_n\), and 
	\(\mathcal H_{cs}^{(0)}=\mathcal H_{cs,n}\).
	
	For a given phase-field variable \(d^{(m)}\), the displacement field is obtained from the weak form of mechanical equilibrium:
	\begin{equation}
		\label{eq:weak_form_u_numerical}
		\begin{aligned}
			R_{\mathbf u}
			&\left(
			\mathbf u^{(m+1)},\delta\mathbf u;d^{(m)}
			\right) \\
			&=
			\int_{\Omega}
			\left[
			g(d^{(m)})\bm{\sigma}_0^+\left(\bm{\varepsilon}(\mathbf u^{(m+1)})\right) + \bm{\sigma}_0^-\left(\bm{\varepsilon}(\mathbf u^{(m+1)})\right)
			\right]
			:
			\bm{\varepsilon}(\delta\mathbf u)\,d\Omega
			-
			\int_{\partial_N\Omega}
			\bar{\mathbf t}\cdot\delta\mathbf u\,dS.
		\end{aligned}
	\end{equation}
	Here, \(\bar{\mathbf t}\) denotes the traction prescribed on the Neumann boundary \(\partial_N\Omega\); this term is active in the traction-controlled example of \cref{subsec:sent_plate} and vanishes in the displacement-controlled tests.
	
	After the displacement solve, the local threshold \(\mathcal H_{cs}^{(m+1)}\) is updated according to the freezing rule in \cref{eq:Hcs_freezing_thermo}. 
	For effective stress states with \(\sigma_1\leq0\), a sufficiently large finite value is used in place of the formal infinite threshold to avoid numerical overflow and suppress tensile crack nucleation in compressive regions.
	
	For fixed \(\mathbf u^{(m+1)}\) and \(\mathcal H_{cs}^{(m+1)}\), the phase-field subproblem is solved in the admissible space 
	\(\mathcal D_{n+1}=\{d\in H^1(\Omega)\mid d_n\leq d\leq1\}\). 
	The weak residual is
	\begin{equation}
		\label{eq:weak_form_d_numerical}
		\begin{aligned}
			R_d &\left( d, \delta d; \mathbf u^{(m+1)}, \mathcal H_{cs}^{(m+1)} \right) \\
			&= \int_{\Omega} \left[ -2(1-d) \left( \psi_e^+\left(\bm{\varepsilon}(\mathbf u^{(m+1)})\right) - \mathcal H_{cs}^{(m+1)} + \psi_c \right) \delta d \right] d\Omega \\
			&\quad + \int_{\Omega} \left[ \frac{3G_c}{8l_0}\delta d + \frac{3G_c l_0}{4} \nabla d\cdot\nabla(\delta d) \right] d\Omega .
		\end{aligned}
	\end{equation}
	
	To improve the robustness of the staggered iteration, the phase-field update \(\widetilde d^{(m+1)}\) is under-relaxed \citep{radtke2016convergence}:
	\begin{equation}
		\label{eq:phase_field_under_relaxation}
		d^{(m+1)}
		=
		(1-\alpha)d^{(m)}
		+
		\alpha\widetilde d^{(m+1)},
		\qquad
		\alpha\in(0,1] .
	\end{equation}
	In the present simulations, \(\alpha\) is chosen empirically in the range \(0.1\)--\(0.5\).
	
	The staggered iteration is accepted when the out-of-balance residuals of both subproblems, evaluated with the most recently updated fields, satisfy
	\begin{equation}
		\label{eq:staggered_convergence}
		\max\Bigl\{
		\bigl\|R_{\mathbf u}\bigl(\mathbf u^{(m+1)},\delta\mathbf u;d^{(m+1)}\bigr)\bigr\|_{L_2},\,
		\bigl\|R_{d}\bigl(d^{(m+1)},\delta d;\mathbf u^{(m+1)},\mathcal H_{cs}^{(m+1)}\bigr)\bigr\|_{L_2}
		\Bigr\} < 10^{-8},
	\end{equation}
	where the phase-field residual is evaluated in the bound-constrained sense, that is, the components associated with the active constraints \(d=d_n\) or \(d=1\) are excluded.
	After convergence, the fields \((\mathbf u_{n+1},d_{n+1})\) are committed, and \(\mathcal H_{cs,n+1}\) is updated according to \cref{eq:Hcs_freezing_thermo}.
	
	All of the numerical examples are carried out in FEniCSx \citep{Baratta2023}. The weak forms \cref{eq:weak_form_u_numerical,eq:weak_form_d_numerical} are written in Unified Form Language \citep{alnaes2014}, and the corresponding consistent Jacobians are obtained by automatic differentiation. The constraints \(d\geq d_n\) and \(d\leq1\) are imposed directly through a bound-constrained solver in PETSc \citep{balay2025petsc}. The complete staggered procedure is summarized in \cref{alg:numerical_scheme}.
	
	\begin{algorithm}[H]
		\caption{Staggered solution scheme with under-relaxation}
		\label{alg:numerical_scheme}
		\begin{algorithmic}[1]
			\State \textbf{Given:} \((\mathbf u_n,d_n,\mathcal H_{cs,n})\) at time step \(t_n\).
			\State \textbf{Initialize:} \(\mathbf u^{(0)}=\mathbf u_n\), \(d^{(0)}=d_n\), \(\mathcal H_{cs}^{(0)}=\mathcal H_{cs,n}\), and \(m=0\).
			\Repeat
			\State Solve the displacement weak form \cref{eq:weak_form_u_numerical} with fixed \(d^{(m)}\) to obtain \(\mathbf u^{(m+1)}\).
			\State Update \(\mathcal H_{cs}^{(m+1)}\) according to the freezing rule \cref{eq:Hcs_freezing_thermo}.
			\State Solve the bound-constrained phase-field problem \cref{eq:weak_form_d_numerical} to obtain \(\widetilde d^{(m+1)}\).
			\State Apply the under-relaxation update in \cref{eq:phase_field_under_relaxation}.
			\State Check the convergence criterion in \cref{eq:staggered_convergence}.
			\State Set \(m\gets m+1\).
			\Until{the convergence criterion is satisfied}
			\State Commit the converged fields \((\mathbf u_{n+1},d_{n+1})\) and update \(\mathcal H_{cs,n+1}\).
		\end{algorithmic}
	\end{algorithm}

	\section{Numerical examples}
	\label{sec:examples}
	
	The present work proposes a shifted energy barrier phase-field formulation for tensile-dominated brittle fracture. 
	The purpose is to introduce the tensile strength as an independent material parameter, reduce the sensitivity of crack nucleation to the regularization length, avoid non-physical damage initiation under compression, and numerically recover the toughness-controlled trend of linear elastic fracture mechanics (LEFM) for sufficiently large cracks. 
	Based on these aims, the numerical examples are arranged as follows.
	
	\begin{itemize}
		\item The uniaxial tension test of a slender bar, presented in \cref{subsec:uniaxial_tension}, focuses on the decoupling between the prescribed tensile strength \(\sigma_t\) and the regularization length \(l_0\), and verifies the finite element implementation against the closed-form solutions of \cref{sec:1d_analysis}.
		
		\item The two compression examples, presented in \cref{subsec:uniaxial_compression,subsec:hole_compression}, address crack nucleation under compressive loading. 
		The rectangular plate represents a nearly pure compressive stress state, whereas the plate with a central hole introduces local tensile stress concentrations within an overall compressive loading condition.
		
		\item The single edge cracked plate, presented in \cref{subsec:sent_plate}, is considered to study the size effect caused by pre-existing cracks. 
		It provides a check of the transition from strength-controlled failure for short cracks to the Griffith-type toughness-controlled limit, represented here by the LEFM solution, for long cracks.
		
		\item The square plate and the cruciform specimen, presented in \cref{subsec:multiaxial,subsec:cruciform_2d}, are used for multiaxial fracture. 
		The former gives controlled proportional stress paths under plane stress, while the latter introduces a more structural stress distribution under equibiaxial tension.
	\end{itemize}

	\subsection{Uniaxial tension of a bar}
	\label{subsec:uniaxial_tension}
	
	The first example is the uniaxial tension of a bar, a standard benchmark for the nucleation behavior of phase-field models \citep{Wu2017,Feng2021,Zolesi2024,Greco2025}. The bar has length $L = 10$~mm and unit cross-section, and is solved with a one-dimensional finite element model of 200 linear elements ($h = 0.05$~mm). The left end is fixed, $u(0)=0$, and a monotonically increasing displacement $u(L)=\bar{u}(t)$ is applied at the right end, ramped as $\bar{u}(t) = \dot{\bar u}\,t$ with $\dot{\bar u} = 1$~mm/s; the time step is $\Delta t = 1 \times 10^{-6}$~s. The phase field is set to $d=0$ at both ends. The material parameters are $E = 7.0 \times 10^4$~MPa and $G_c = 0.008$~N/mm, representative of a brittle glass. In the absence of stress concentrations, failure is governed entirely by nucleation, and the computed response serves to verify the finite element implementation against the closed-form solutions of \cref{sec:1d_analysis}.
	
	\begin{figure}[htbp]
		\centering
		\includegraphics[width=1.0\textwidth]{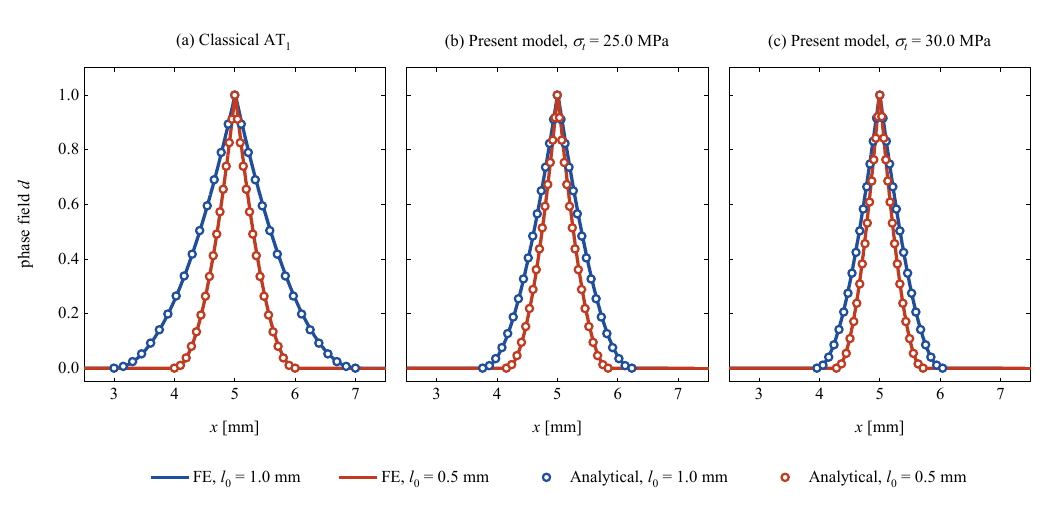}
		\caption{Ultimate phase-field profiles in the uniaxial tension test
			($l_0 = 1.0$ and $0.5$~mm): finite element solutions (solid lines)
			and closed-form solutions of \cref{sec:1d_analysis} (open circles),
			for (a) the classical $\mathrm{AT}_1$ model and (b, c) the present
			model at $\sigma_t = 25$ and $30$~MPa.}
		\label{fig:1d_profiles_validation}
	\end{figure}
	
	Six simulations are performed: the present model with $\sigma_t = 25$~MPa and
	$\sigma_t = 30$~MPa, and the classical $\mathrm{AT}_1$ model, each with
	$l_0 = 1.0$~mm and $l_0 = 0.5$~mm. All four cases of the present model satisfy
	the admissibility condition \cref{eq:l0_condition_thermo}. The analytical
	characterization of the six cases is summarized in \cref{tab:1d_summary}: the
	fracture stress, the shift parameter $\gamma$ of \cref{eq:gamma_definition},
	which ranges from $0.488$ to $3.286$, and the half-width $D$ of the
	localization band given by \cref{eq:1d_bandwidth}.

	\begin{table}[htbp]
		\centering
		\caption{Analytical characterization of the one-dimensional localized
			solutions for the classical $\mathrm{AT}_1$ model and the present
			model ($E = 7.0 \times 10^4$~MPa, $G_c = 0.008$~N/mm). The fracture stress
			$\sigma_f$ equals $\sigma_c = \sqrt{3 E G_c / (8 l_0)}$ for
			$\mathrm{AT}_1$ and the prescribed strength $\sigma_t$ for the
			present model; $D$ denotes the half-width of the localization band,
			\cref{eq:1d_bandwidth}.}
		\label{tab:1d_summary}
		\begin{tabular}{llcccc}
			\toprule
			Model & $l_0$ [mm] & $\sigma_f$ [MPa] & $\gamma$
			& $D$ [mm] & $D/l_0$ \\
			\midrule
			$\mathrm{AT}_1$ & 1.0 & 14.49 & 0     & 2.000 & 2.00 \\
			$\mathrm{AT}_1$ & 0.5 & 20.49 & 0     & 1.000 & 2.00 \\
			Present         & 1.0 & 25.00 & 1.976 & 1.236 & 1.24 \\
			Present         & 0.5 & 25.00 & 0.488 & 0.844 & 1.69 \\
			Present         & 1.0 & 30.00 & 3.286 & 1.042 & 1.04 \\
			Present         & 0.5 & 30.00 & 1.143 & 0.718 & 1.44 \\
			\bottomrule
		\end{tabular}
	\end{table}

	\Cref{fig:1d_profiles_validation} compares the ultimate phase-field profiles
	along the bar with the closed-form solutions of \cref{sec:1d_analysis}. The
	profiles computed with the classical $\mathrm{AT}_1$ model reproduce the
	parabolic solution \cref{eq:at1_profile} with support half-width $2 l_0$,
	while those of the present model follow the cosine-type profile
	\cref{eq:1d_profile}, with the support agreeing with the analytical half-width
	$D$ of \cref{eq:1d_bandwidth} in all six cases. The results confirm the two
	features of the localized solution established in \cref{sec:1d_analysis}:
	at fixed $l_0$, the band narrows as the prescribed strength increases
	($D = 1.24$~mm and $1.04$~mm for $\sigma_t = 25$ and $30$~MPa at
	$l_0 = 1.0$~mm); and at fixed strength, the band width scales sub-linearly
	with $l_0$ ($D/l_0 = 1.24$ at $l_0 = 1.0$~mm versus $1.69$ at $l_0 = 0.5$~mm
	for $\sigma_t = 25$~MPa), in contrast with the proportional scaling
	$D = 2 l_0$ of the $\mathrm{AT}_1$ model.

\begin{figure}[htbp]
		\centering
		\includegraphics[width=1.0\textwidth]{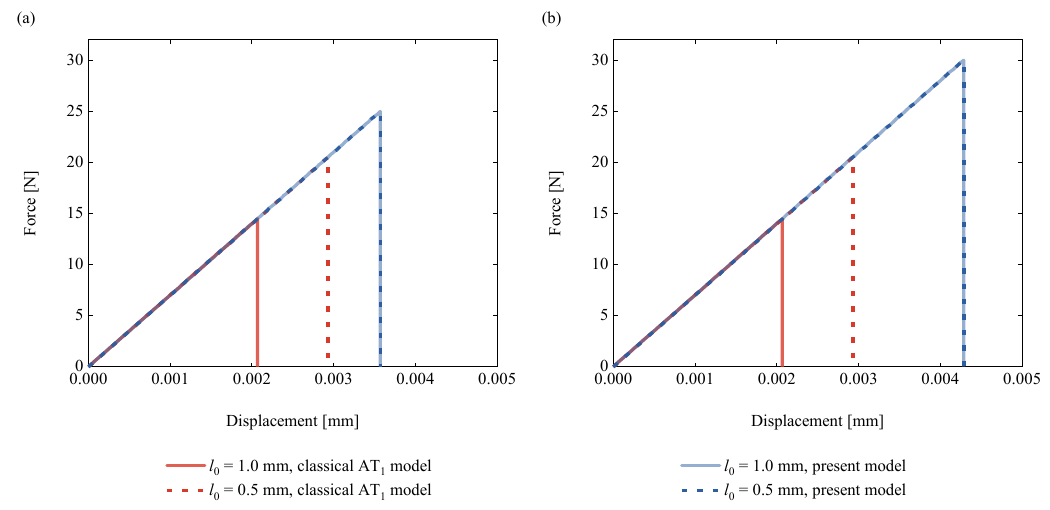}
		\caption{Force--displacement response of the bar under uniaxial tension
			for the present model (blue) and the classical $\mathrm{AT}_1$ model
			(red), with $l_0 = 1.0$~mm (solid) and $l_0 = 0.5$~mm (dashed):
			(a)~$\sigma_t = 25$~MPa; (b)~$\sigma_t = 30$~MPa. The two blue curves
			coincide in each panel, with peak load set by the prescribed strength;
			the red curves peak at $\sigma_c = \sqrt{3 G_c E / (8 l_0)}$,
			identical in both panels.}
		\label{fig:uniaxial_response}
	\end{figure}

	\Cref{fig:uniaxial_response} shows the computed force--displacement responses.
	All curves follow the same elastic branch and fail abruptly at their
	respective peak loads. For the classical $\mathrm{AT}_1$ model, the peak
	stress follows the intrinsic value $\sigma_c$ of \cref{rem:inherent_barrier}:
	it increases from $14.49$~MPa to $20.49$~MPa when $l_0$ is halved, and is
	identical in both panels regardless of the intended strength. For the present
	model, the two regularization lengths yield coinciding curves whose peak
	stress equals the prescribed strength. This insensitivity to $l_0$ follows
	directly from the construction of the shifted driving force: by
	\cref{eq:kkt_dual}, damage initiates when $\psi_e^+$ reaches the frozen
	threshold $\mathcal{H}_{cs}$ of \cref{eq:Hcs_freezing_thermo}, which in the
	present uniaxial setting equals $\sigma_t^2/(2E)$ by \cref{eq:1d_threshold},
	so that the onset condition reduces to $\sigma = \sigma_t$ irrespective of
	$l_0$; the regularization length enters the onset condition only through the
	admissibility bound \cref{eq:l0_condition_thermo}. In the classical
	$\mathrm{AT}_1$ model, by contrast, the threshold is the intrinsic barrier
	$\psi_c = 3 G_c / (16 l_0)$ of \cref{eq:inherent_barrier}, which ties the
	nucleation stress to the regularization length. The comparison thus isolates
	the central property of the present formulation: the nucleation stress is a
	material input rather than a byproduct of $l_0$.

	\subsection{Uniaxial compression of a rectangular plate}
	\label{subsec:uniaxial_compression}
	A standard uniaxial compression benchmark \citep{Feng2021,li2026bridging} is simulated to examine the model's response under macroscopic compressive loading. While classical phase-field models typically capture shear-driven damage under such conditions, this benchmark specifically verifies the present formulation's capacity to isolate pure tension-driven failure.
	
\begin{figure}[htbp]
		\centering
		\includegraphics[width=1.0\textwidth]{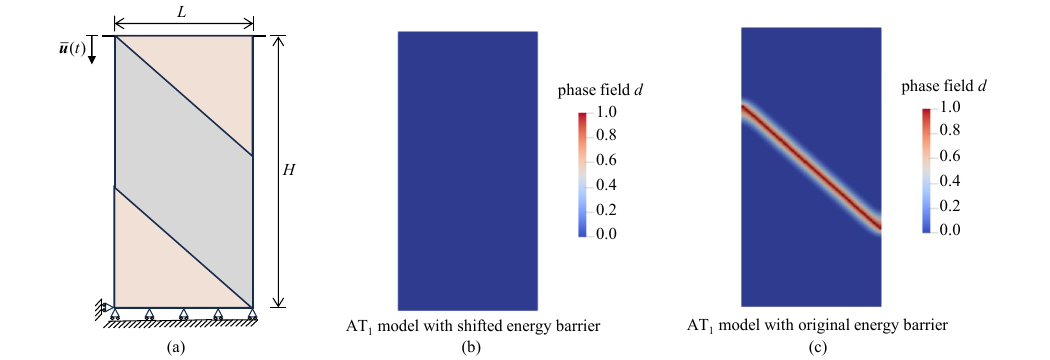}
		\caption{Phase-field predictions for the uniaxial compression test. (a) Geometry and boundary conditions. Final phase-field contours illustrate the differing damage mechanisms: (b) the shifted energy barrier formulation isolates pure tension-driven failure, remaining intact under compression, whereas (c) the classical \(\mathrm{AT}_1\) model captures a diagonal shear band.}
		\label{fig:compression_test}
	\end{figure}

	The geometric configuration and boundary conditions are illustrated in \cref{fig:compression_test}(a). The rectangular specimen has a width $L = 10.0$~mm and a height $H = 20.0$~mm. The bottom edge is constrained vertically, with the bottom-left corner fixed to eliminate rigid body translation. A monotonically increasing compressive displacement is applied to the top boundary, ramped as $\bar{u}(t) = \dot{\bar u}\,t$ with $\dot{\bar u} = 1$~mm/s. The time step used in the simulation is $\Delta t = 1 \times 10^{-3}$~s. To prevent damage initiation from the stress concentrations at the boundaries, a Dirichlet condition $d=0$ is enforced within the beige-shaded regions. 
	
	The material properties are: Young's modulus $E = 7.0 \times 10^4$~MPa, Poisson's ratio $\nu = 0.20$, critical energy release rate $G_c = 0.008$~N/mm, and tensile strength $\sigma_t = 25.0$~MPa. The regularization length is $l_0 = 0.5$~mm. The domain is discretized using a uniform rectangular mesh with a characteristic size of $h = 0.05$~mm, providing a spatial resolution of $h = l_0 / 10$.
	
	The final phase-field distributions are compared in \cref{fig:compression_test}(b) and (c). For an idealized pure tension-driven brittle material, damage initiation is not expected under unconfined uniaxial compression. As depicted in \cref{fig:compression_test}(b), the proposed formulation maintains a zero-damage state throughout the loading history. In contrast, the classical \(\mathrm{AT}_1\) model based on standard strain energy decomposition yields a diagonal shear band (\cref{fig:compression_test}(c)). This comparison demonstrates that the shifted energy barrier mechanism explicitly separates the tension-driven failure mode from compression-induced damage evolution.

	\subsection{Uniaxial compression of a plate with central hole}
	\label{subsec:hole_compression}

	A rectangular plate with a central hole under uniaxial compression \citep{li2026bridging,feng2023unified} is simulated to evaluate the model's response to structural stress concentrations. This benchmark examines crack nucleation within a non-uniform stress field without pre-existing flaws. 
	
\begin{figure}[htbp]
		\centering
		\includegraphics[width=1\textwidth]{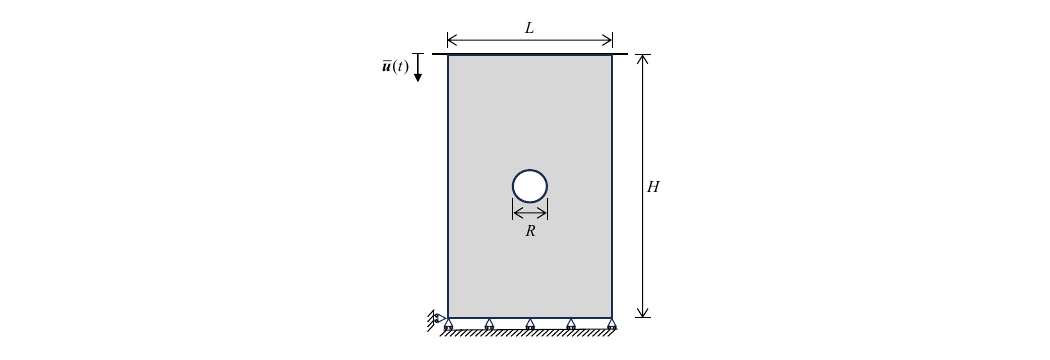}
		\caption{Schematic of the uniaxial compression test for a rectangular plate with a central hole: geometry and boundary conditions.}
		\label{fig:hole_setup}
	\end{figure}

	The geometric configuration and boundary conditions are illustrated in \cref{fig:hole_setup}. The plate dimensions are width $L = 100.0$~mm and height $H = 170.0$~mm, with a central circular hole of radius $R = 7.5$~mm. The bottom edge is constrained vertically, and the bottom-center node is fixed to eliminate rigid body translation. A monotonically increasing compressive displacement is applied to the top edge, ramped as $\bar{u}(t) = \dot{\bar u}\,t$ with $\dot{\bar u} = 1$~mm/s. The time step used in the simulation is $\Delta t = 1 \times 10^{-4}$~s. The lateral edges and the hole surface remain traction-free.
	
	The material properties are: Young's modulus $E = 2.1 \times 10^5$~MPa, Poisson's ratio $\nu = 0.3$, critical energy release rate $G_c = 0.001$~N/mm, and tensile strength $\sigma_t = 12.0$~MPa. The regularization length is $l_0 = 2.0$~mm. The two-dimensional computational domain is discretized using quadrilateral elements. The mesh in the anticipated crack propagation regions is refined to $h = 0.2$~mm, ensuring a spatial resolution of $h = l_0 / 10$.

\begin{figure}[htbp]
		\centering
		\includegraphics[width=\textwidth]{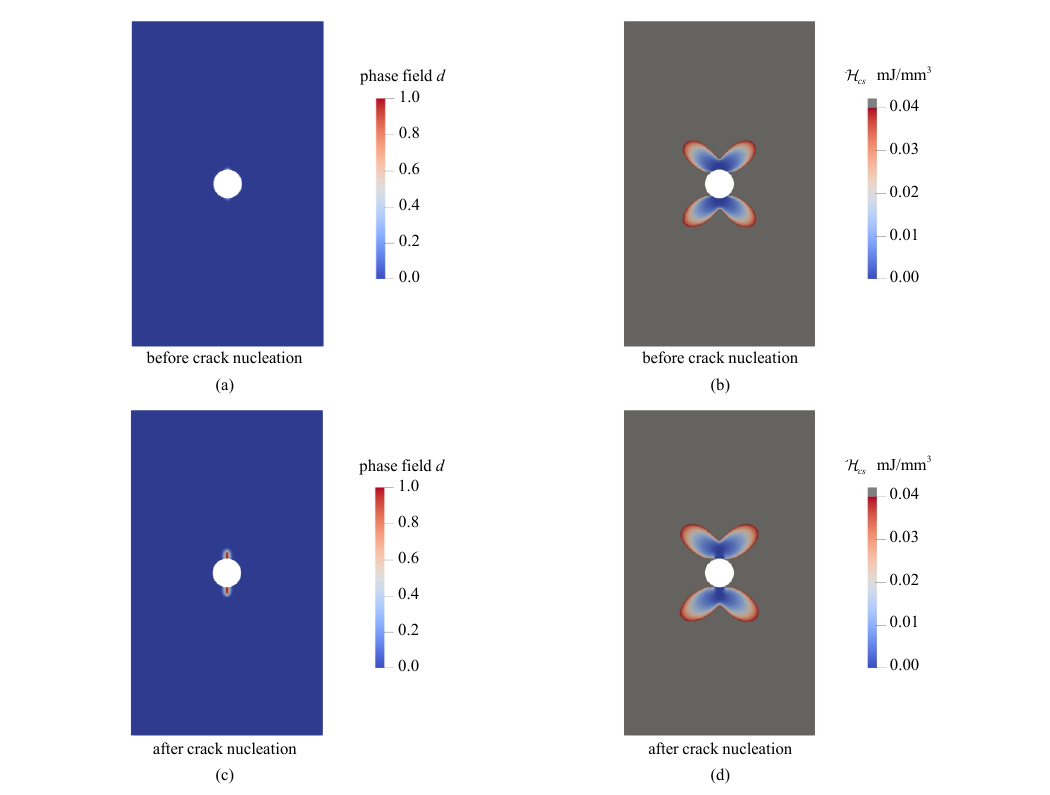} 
		\caption{Phase-field damage and energy threshold evolution under uniaxial compression. Contours of the phase-field variable $d$ and the frozen energy threshold $\mathcal{H}_{cs}$ are shown immediately before crack nucleation (a, b) and after the initiation of axial splitting cracks at the poles (c, d). The elevated $\mathcal{H}_{cs}$ values at the compressed lateral equators in (b) and (d) suppress crack nucleation in these regions.}
		\label{fig:hole_contours}
	\end{figure}

	The evolution of the phase-field damage and the energy threshold during loading
	is presented in \cref{fig:hole_contours}. As shown in
	\cref{fig:hole_contours}(a) and (c), the macroscopic compressive load
	induces local tensile stress concentrations at the top and bottom poles of the
	hole, which drive the nucleation and propagation of axial splitting cracks.
	
	The mechanism by which the compressed regions remain intact is a distinctive feature of the
	present formulation. In classical phase-field models, the lateral equators of
	the hole accumulate a large compressive strain energy that is partly retained in
	the active energy through the tension--compression split, and can therefore
	trigger spurious damage in these strongly compressed zones. Here, the resistance
	to nucleation is instead governed by the spatially varying,
	stress-state-dependent threshold $\mathcal{H}_{cs}$. Because $\mathcal{H}_{cs}$ scales the
	active energy to the Rankine surface through the factor $(\sigma_t/\sigma_1)^2$
	in \cref{eq:psi_cs}, it takes small values where a tensile concentration
	develops, namely at the poles, and large values where the maximum effective
	principal stress is low or non-positive, namely at the compressed equators,
	diverging for $\sigma_1 \le 0$ (see \cref{rem:compressive_threshold}). This
	spatial pattern is shown in \cref{fig:hole_contours}(b) and (d), where
	$\mathcal{H}_{cs}$ is markedly elevated along the equators and low at the poles.
	As a result, the shifted driving force $\psi_e^+ - \mathcal{H}_{cs} + \psi_c$ can
	vanish only at the poles, whereas at the equators it remains negative so that
	damage cannot initiate. The elevated threshold thus acts as an intrinsic,
	stress-driven gate that confines crack nucleation to regions of active tensile
	stress, ensuring that structural failure is driven exclusively by local tensile
	stresses rather than by compressive strain energy.
	
	\subsection{Single edge cracked plate under tension}
	\label{subsec:sent_plate}
	
	A size effect study of a single edge cracked plate under uniform tension \citep{chockalingam2026phase,sargado2018high,kristensen2021assessment} is conducted to examine the transition from strength-governed to toughness-governed failure as a function of the initial crack length.
	
	\begin{figure}[htbp]
		\centering
		\includegraphics[width=\textwidth]{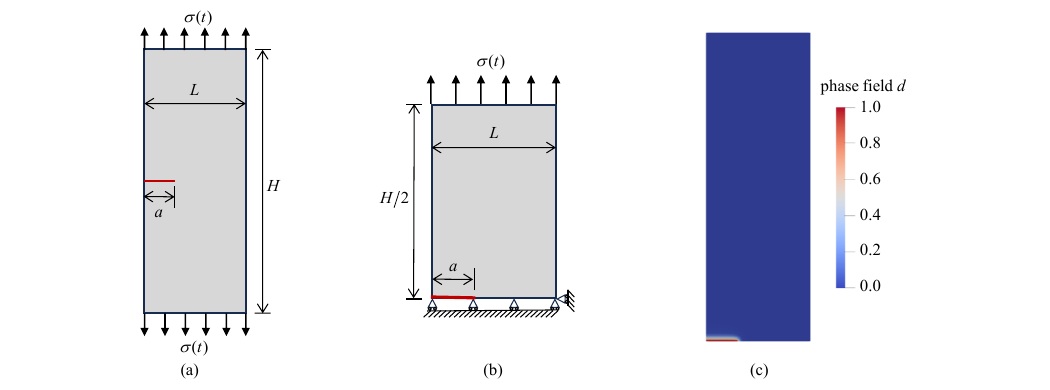} 
		\caption{Single edge cracked plate under tension. (a) Full geometry and loading conditions for a plate with an initial crack of length $a$ subjected to a uniform boundary traction $\sigma(t)$. (b) Equivalent half-symmetry computational model and boundary conditions. (c) Regularized phase-field representation of the initial crack using the damage variable $d$.}
		\label{fig:sent_setup}
	\end{figure}
	
	The geometric configuration is detailed in \cref{fig:sent_setup}(a). To reduce computational cost, the half-symmetry model shown in \cref{fig:sent_setup}(b) is employed. The specimen dimensions are width $L = 100.0$~mm and half-height $H/2 = 300.0$~mm, with a horizontal initial crack of length $a$. Symmetry constraints ($u_y = 0$) are enforced along the uncracked ligament of the bottom edge. A monotonically increasing uniform tensile traction $\sigma(t) = \dot\sigma\,t$ with $\dot\sigma = 1$~MPa/s is applied to the top boundary under plane strain conditions. The time step used in the simulation is $\Delta t = 1 \times 10^{-3}$~s. The material properties are: Young's modulus $E = 3.0 \times 10^4$~MPa, Poisson's ratio $\nu = 0.23$, critical energy release rate $G_c = 0.03$~N/mm, and tensile strength $\sigma_t = 17.5$~MPa. Simulations are performed for two regularization lengths, $l_0 = 4.0$~mm and $l_0 = 5.0$~mm. The domain is discretized using quadrilateral elements with a uniform mesh size of $h = l_0/10$.
	
	For this geometry, the critical failure stress $\sigma_f$ is bounded by two limiting mechanisms. In the strength-governed limit for short cracks ($a/L \to 0$), the failure stress approaches the material strength:
	\begin{equation}
		\sigma_f = \sigma_t.
		\label{eq:strength_limit}
	\end{equation}
	In the toughness-governed limit for long cracks, the response follows the LEFM solution \citep{anderson2005fracture}. Under the plane strain condition, the critical stress is:
	\begin{equation}
		\sigma_f = \frac{1}{F(a/L)} \sqrt{\frac{E G_c}{\pi (1-\nu^2)a}},
		\label{eq:lefm_limit}
	\end{equation}
	where $F(a/L)$ is the dimensionless geometric correction function \citep{Tada2000}:
	\begin{equation}
		F(a/L) = \frac{1}{\cos\left(\frac{\pi a}{2L}\right)} \left[ 0.752 + 2.02 \frac{a}{L} + 0.37 \left(1 - \sin\frac{\pi a}{2L}\right)^3 \right] \sqrt{\frac{2L}{\pi a} \tan\left(\frac{\pi a}{2L}\right)}.
	\end{equation}
	
	Traction control is used here because only the critical stress $\sigma_f$ is sought. Under the prescribed traction $\sigma(t)$, the pre-peak response is stable and converges at every increment, whereas the onset of macroscopic fracture triggers a structural instability that manifests as a loss of convergence of the quasi-static solution scheme. The failure stress $\sigma_f$ is therefore extracted as the traction magnitude at the last converged load step.
	
	\begin{figure}[htbp]
		\centering
		\includegraphics[width=1.0\textwidth]{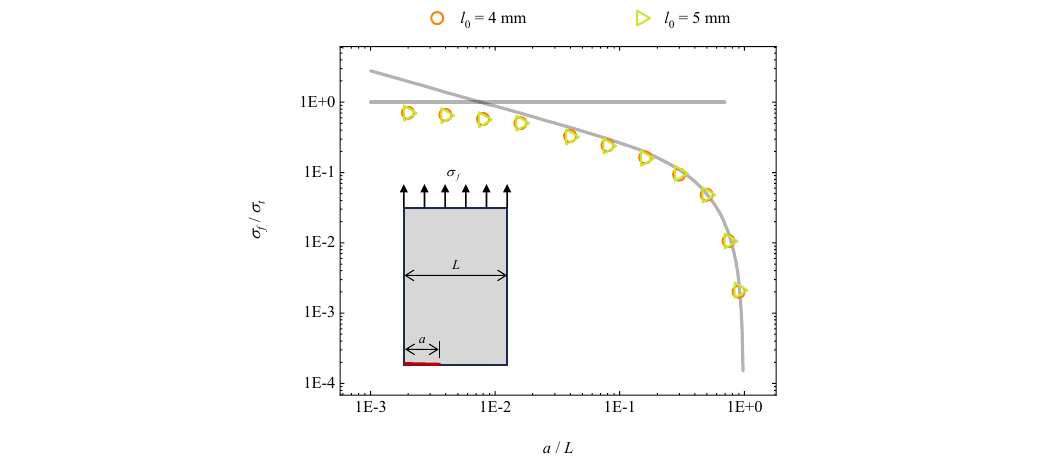} 
		\caption{Size effect and transition of failure mechanisms in the single edge cracked plate. Normalized failure stress $\sigma_f / \sigma_t$ is plotted against the dimensionless initial crack length $a/L$ on a logarithmic scale. Numerical results for $l_0 = 4.0$~mm and $l_0 = 5.0$~mm are compared with the theoretical strength limit (horizontal gray line) and the LEFM toughness limit (sloped gray line).}
		\label{fig:size_effect}
	\end{figure}
	
	\Cref{fig:size_effect} plots the computed size effect curves alongside the theoretical limits. For small $a/L$, the predicted failure stress converges to the strength limit \cref{eq:strength_limit} defined by the Rankine criterion. As $a/L$ increases, the failure stress transitions to the toughness-dominated regime, aligning with the LEFM asymptotic solution \cref{eq:lefm_limit}. The predictions for $l_0 = 4.0$~mm and $l_0 = 5.0$~mm are in close agreement across the investigated crack length range. The shifted energy barrier formulation thus decouples the macroscopic failure load from the regularization length, ensuring a consistent transition between the strength and Griffith fracture limits.
	
	\subsection{Multiaxial fracture envelope of a square plate under plane stress}
	\label{subsec:multiaxial}
	
\begin{figure}[htbp]
		\centering
		\includegraphics[width=1.0\textwidth]{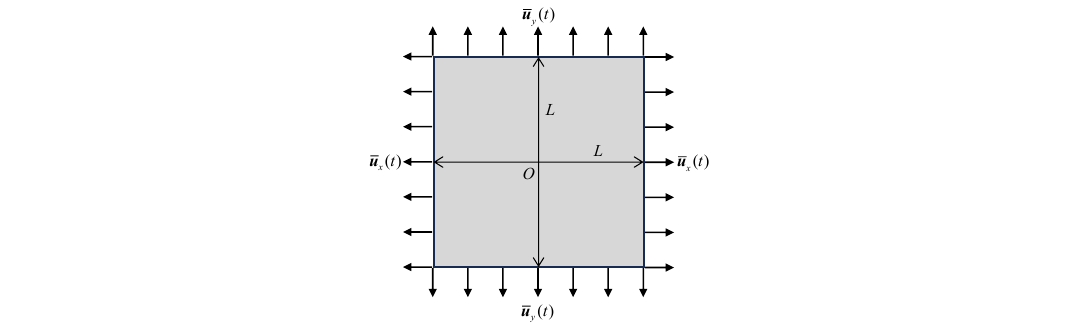} 
		\caption{Schematic of the multiaxial tension test for a square plate: geometry and boundary conditions. The domain is subjected to proportional outward displacements on all four edges, governed by the loading angle $\Theta$.}
		\label{fig:equibiaxial_setup}
	\end{figure}

	A square plate under proportional loading \citep{Vicentini2024,chockalingam2026phase} is simulated to evaluate the model's performance in multiaxial stress states. The geometry and boundary conditions are detailed in \cref{fig:equibiaxial_setup}. The plate has a side length $L = 5.0$~mm and is modeled under the plane stress assumption. Plane stress is adopted so that the out-of-plane stress vanishes ($\sigma_3=0$) and the prescribed in-plane principal stresses map directly onto the Rankine envelope in the $(\sigma_1,\sigma_2)$ plane. Material parameters are: Young's modulus $E = 100.0$~MPa, Poisson's ratio $\nu = 0.3$, critical energy release rate $G_c = 0.2$~N/mm, and tensile strength $\sigma_t = 8.0$~MPa. Two regularization lengths, $l_0 = 0.4$~mm and $l_0 = 0.2$~mm, are investigated. The domain is discretized using quadrilateral elements with a characteristic size $h = l_0 / 4$. 
	
\begin{figure}[htbp]
		\centering
		\includegraphics[width=1.0\textwidth]{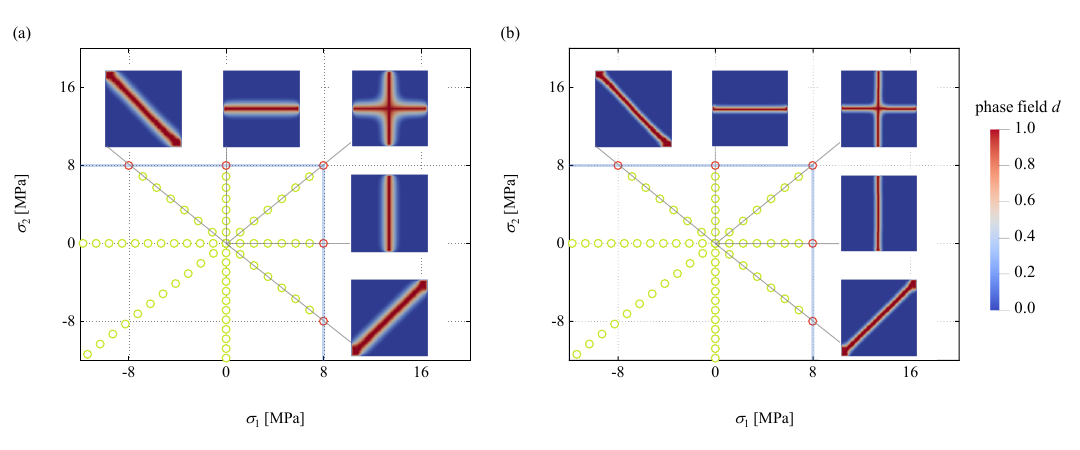}
		\caption{
			Macroscopic failure envelope and corresponding phase-field fracture patterns in the principal stress space ($\sigma_1$, $\sigma_2$). The numerical crack nucleation points (red circles) under various proportional loading paths (light green circles) are compared with the theoretical Rankine strength limit ($\sigma_t = 8$~MPa, blue solid lines). The simulation results for two different regularization lengths, (a) $l_0 = 0.4$~mm and (b) $l_0 = 0.2$~mm, overlap at the theoretical boundary. The inset contours show that the model captures distinct crack topologies, such as pure tensile splitting and equibiaxial orthogonal cracking, insensitive to $l_0$ in the tested admissible range.
		}
		\label{fig:failure_envelope}
	\end{figure}

	Proportional displacements are applied to the four outer edges. The horizontal and vertical displacement components are mathematically defined as:
	\begin{equation}
		\bar{u}_x(\pm L/2, y) = \pm \frac{1}{2}\bar{u}(t) \cos\Theta, \quad \bar{u}_y(x, \pm L/2) = \pm \frac{1}{2}\bar{u}(t) \sin\Theta,
		\label{eq:proportional_loading}
	\end{equation}
	where $\bar{u}(t)$ is the monotonically increasing generalized displacement, and $\Theta$ denotes the loading angle. The generalized displacement is ramped as $\bar{u}(t) = \dot{\bar u}\,t$ with $\dot{\bar u} = 1$~mm/s, and the time step used in the simulation is $\Delta t = 1 \times 10^{-3}$~s. To prevent damage nucleation from boundary stress concentrations, a Dirichlet condition $d=0$ is enforced along the entire perimeter. Consequently, crack propagation arrests near the edges. Varying the displacement angle $\Theta$ generates different linear loading paths in the principal stress space ($\sigma_1$, $\sigma_2$), indicated by the light green circles in \cref{fig:failure_envelope}. Due to the Poisson effect under plane stress, the resulting principal stress ratio differs from the prescribed displacement ratio, following the relation $\sigma_2/\sigma_1 = (\sin\Theta + \nu \cos\Theta)/(\cos\Theta + \nu \sin\Theta)$.
	
	\Cref{fig:failure_envelope} plots the computed failure envelope alongside the phase-field crack patterns. Across all loading paths, the numerical crack nucleation points (red circles) align with the theoretical Rankine strength limit (solid blue lines). The predictions for $l_0 = 0.4$~mm and $l_0 = 0.2$~mm nearly coincide, indicating the length-scale insensitivity of the crack initiation limit. The model naturally captures distinct crack topologies dictated by the local stress state. Uniaxial tension produces a single straight splitting crack perpendicular to the principal loading direction. Equibiaxial tension ($\sigma_1 = \sigma_2$) yields symmetric orthogonal cross-shaped cracks. Mixed tension--compression states (second and fourth quadrants) generate inclined cracks, slightly broader than those under pure tension. Under compressive loading paths (third quadrant), no cracks nucleate. Overall, these results demonstrate that the proposed framework enforces the Rankine fracture criterion in a manner insensitive to the phase-field regularization length.

	\subsection{Equibiaxial tension of a cruciform specimen}
	\label{subsec:cruciform_2d}
	
\begin{figure}[htbp]
		\centering
		\includegraphics[width=1.0\textwidth]{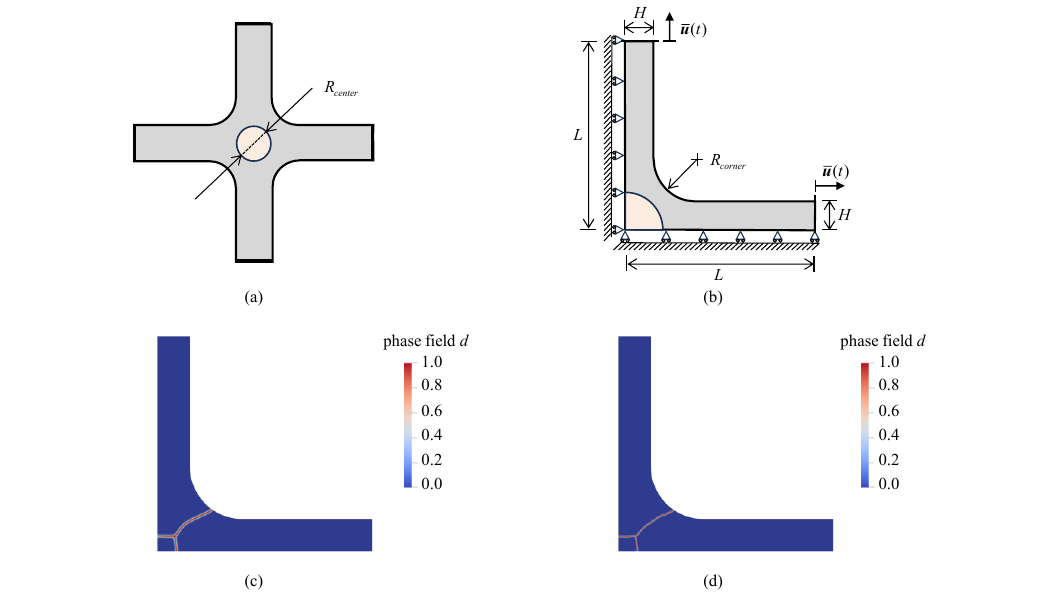} 
		\caption{Numerical setup and phase-field predictions for the cruciform specimen under equibiaxial tension. (a) Schematic of the full geometry highlighting the central thinned region of radius $R_{\mathrm{center}}$. (b) Quarter-symmetry computational domain detailing the geometric dimensions, symmetric boundary constraints, and prescribed uniform outward displacements $\bar{u}(t)$. Final phase-field damage contours illustrating the predicted crack trajectories for regularization lengths of (c) $l_0 = 1.0$~mm and (d) $l_0 = 0.5$~mm.}
		\label{fig:cruciform}
	\end{figure}

	We simulate a cruciform specimen subjected to equibiaxial tension, adopting a plane stress formulation and a quarter-symmetry domain. The geometry is defined by an arm length $L = 165.0$~mm, an arm half-width $w = 30.0$~mm, a fillet radius $R_{\mathrm{corner}} = 40.0$~mm, and a central region radius $R_{\mathrm{center}} = 25.0$~mm, as illustrated in \cref{fig:cruciform}(a) and (b). 
	
	In experimental mechanics, the central gage section of a cruciform specimen is routinely machined to a reduced thickness compared to the loading arms. This geometric modification mitigates premature stress concentrations at the re-entrant corners, induces a homogeneous multiaxial stress state, and dictates that macroscopic crack nucleation initiates strictly within the designated central region. To computationally represent this three-dimensional feature within a two-dimensional plane stress framework, a spatially varying thickness function $t_{\mathrm{th}}(x, y)$ is introduced:
	\begin{equation}
		t_{\mathrm{th}}(x, y) = 
		\begin{cases} 
			t_{\mathrm{center}}, & \text{if } x^2 + y^2 \le R_{\mathrm{center}}^2 ,\\ 
			t_{\mathrm{arm}},    & \text{otherwise} ,
		\end{cases}
		\label{eq:thickness}
	\end{equation}
	where $x$ and $y$ denote the Cartesian coordinates relative to the specimen centroid. The thickness is prescribed as $t_{\mathrm{center}} = 1.0$~mm within the central circular region and $t_{\mathrm{arm}} = 5.0$~mm throughout the remainder of the domain. This variable thickness field is then directly incorporated as a scalar multiplier into the volume integrals of the weak forms for both the mechanical equilibrium and the phase-field damage evolution.
	
	The material parameters are set as follows: Young's modulus $E = 100.0$~MPa, Poisson's ratio $\nu = 0.3$, critical energy release rate $G_c = 0.2$~N/mm, and tensile strength $\sigma_t = 6.0$~MPa. To investigate the sensitivity of the formulation to the regularization length, simulations are conducted for $l_0 = 1.0$~mm and $l_0 = 0.5$~mm. The central region is discretized with a characteristic mesh size $h = 0.1$~mm to ensure adequate resolution of the phase-field profile. The prescribed uniform outward displacement is ramped as $\bar{u}(t) = \dot{\bar u}\,t$ with $\dot{\bar u} = 1$~mm/s. The time step used in the simulation is $\Delta t = 2 \times 10^{-2}$~s.
	
\begin{figure}[htbp]
		\centering
		\includegraphics[width=1.0\textwidth]{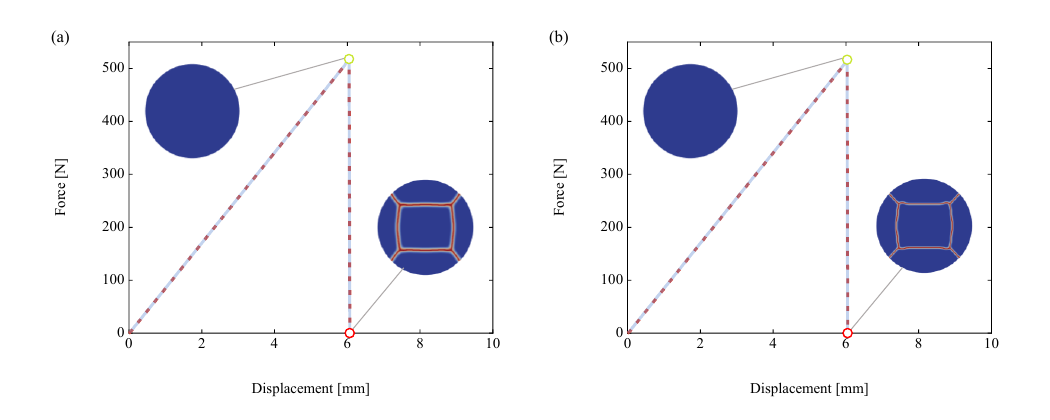}
		\caption{Force--displacement responses of the cruciform specimen under equibiaxial tension for (a) $l_0 = 1.0$~mm and (b) $l_0 = 0.5$~mm. Both cases exhibit a linear elastic regime terminating in an abrupt load drop.}
		\label{fig:cruciform_curve}
	\end{figure}

	The final phase-field fracture patterns and the global force--displacement responses for the two regularization lengths are presented in \cref{fig:cruciform,fig:cruciform_curve}, respectively. \Cref{fig:cruciform}(c) and (d) detail the crack trajectories within the quarter-symmetry computational domain. Crack nucleation initiates at the geometric boundary separating the central thinned region and the thicker loading arms. Subsequently, the damage propagates outward along the horizontal and vertical axes of symmetry. The predicted fracture paths for $l_0 = 1.0$~mm and $l_0 = 0.5$~mm are virtually indistinguishable. The uncracked regions exhibit no spurious damage diffusion, confirming that the shifted energy barrier logic prevents non-physical damage growth prior to ultimate failure. \Cref{fig:cruciform_curve} plots the corresponding macroscopic mechanical behavior. The force--displacement response exhibits a linear elastic regime terminating in an abrupt load drop. The critical failure load and the ultimate displacement are insensitive to the phase-field regularization length. The insets in \cref{fig:cruciform_curve} illustrate the local phase-field state within the central region, reconstructed to the full circular domain via mirror symmetry. The central gage section remains intact during the linear loading phase. Upon reaching the critical displacement, a cross-shaped crack pattern forms, splitting the domain into four symmetric quadrants.

	\section{Conclusions}
	\label{sec:conclusions}
	
	The classical \(\mathrm{AT}_1\) phase-field model couples the crack nucleation stress to the critical energy release rate and the regularization length through an intrinsic energy barrier. To address this artificial coupling, this work presents a shifted energy barrier formulation for tensile-dominated brittle fracture. By mapping the multiaxial Rankine stress criterion to a local active elastic energy threshold and introducing it into the damage driving force, the macroscopic tensile strength is prescribed as an independent material parameter. The shifted threshold, introduced through a restricted variational principle \citep{rosen1954solution,rosen1954use}, renders strength-controlled crack initiation insensitive to $l_0$ within the admissible range, while the standard degraded elastic stress response, the stiffness-degradation structure and the \(\mathrm{AT}_1\) crack surface density are explicitly preserved. A core theoretical feature of the formulation is its conditional thermodynamic admissibility under $\mathcal{H}_{cs} \geq \psi_c$. The Rankine-equivalent threshold functions as a dissipative resistance rather than an additional recoverable Helmholtz energy. By introducing a history-based freezing rule post-initiation, the formulation structurally locks the activated damage resistance, preventing unphysical threshold variations during local unloading or non-proportional loading. Non-negative damage dissipation is ensured under the irreversibility condition $\dot d \ge 0$ in \cref{eq:kkt_primal} together with the admissibility condition $\mathcal{H}_{cs} \ge \psi_c$ in \cref{eq:Hcs_admissibility}. 
	
	In the one-dimensional traction problem the Rankine-equivalent threshold reduces to the constant energy level $\sigma_t^2/(2E)$, and closed-form solutions are derived for the homogeneous softening response and for the localized damage profile. The peak nominal stress equals the prescribed strength, and the localized profile is of cosine type. Its half-width is strictly smaller than the classical value $2l_0$ and grows sub-linearly with the regularization length at large shift, and the parabolic $\mathrm{AT}_1$ profile is recovered as the shift parameter tends to zero. These solutions quantify the effect of the barrier shift and serve as analytical benchmarks for the finite element implementation.
	
	Numerical assessments verify the behavior of the proposed framework across distinct fracture regimes. The main structural and physical characteristics demonstrated by the simulations are summarized as follows:
	
	\begin{itemize}
		\item \textit{Length-scale insensitivity}: Under uniaxial tension, the macroscopic failure load is governed by the prescribed material strength, decoupling the crack nucleation limit from the phase-field regularization length. The computed damage profiles and peak loads reproduce the closed-form solutions, including the cosine-type profile and the analytical half-width of the localization band.
		
		\item \textit{Suppression of unphysical damage in compression}: Under macroscopic compressive loading, the spatial distribution of the local energy threshold suppresses damage nucleation in highly compressed regions, ensuring that fracture localizes exclusively at active tensile stress concentrations.
		
		\item \textit{Transition of failure mechanisms}: For pre-cracked structures, the formulation captures the underlying size effect \citep{kamarei2026nine}. It naturally reproduces the transition from strength-controlled failure for short cracks to the classical toughness-controlled LEFM limit for macroscopic flaws.
		
		\item \textit{Multiaxial crack nucleation}: Under complex loading configurations, including proportional plane-stress paths and non-uniform structural stress fields in a cruciform specimen, the model enforces the Rankine nucleation envelope and yields physically consistent crack initiation patterns.
	\end{itemize}
	
	The present formulation is restricted to tensile-dominated brittle fracture and relies on a length-scale lower bound to ensure non-negative damage dissipation. Extending this approach to compressive-shear failure requires mapping generalized multiaxial strength criteria onto the active energy threshold. Integrating the shifted barrier with directional energy decomposition \citep{feng2023unified} may also provide a potential route to resolve physically consistent crack nucleation orientations under complex stress states. In addition, high-order phase-field approximations can further reduce spatial discretization errors and improve solver efficiency \citep{greco2024higher,greco2026fourth}. To achieve the requisite high-order spatial continuity, isogeometric analysis has emerged as a promising alternative numerical framework \citep{hughes2005isogeometric}, facilitating the application of the proposed model to large-scale engineering problems \citep{morganti2015patient}.
	
	\section*{Acknowledgements}
	
	The authors gratefully acknowledge the financial support provided by the National Natural Science Foundation of China (NSFC) (Grant Nos.\ 12172103, 12572087, and 12020101001), and the Heilongjiang Touyan Innovation Team Program. We also thank Associate Professor Ye Feng at Northwestern Polytechnical University for his helpful discussions on the energy decomposition methods. 
	
	\bibliographystyle{elsarticle-harv}
	\bibliography{references}
	
\end{document}